\documentclass[a4paper,12pt,twoside]{article}
\usepackage[francais]{babel}
\usepackage{amssymb, amsmath, eucal}
\usepackage[all]{xy}
\usepackage{mathrsfs}
\usepackage[dvips]{graphics}
\begin{document}
\textwidth15.5cm
\textheight22.5cm
\voffset=-13mm
\newtheorem{The}{Theorem}[section]
\newtheorem{Lem}[The]{Lemma}
\newtheorem{Prop}[The]{Proposition}
\newtheorem{Cor}[The]{Corollary}
\newtheorem{Rem}[The]{Remark}
\newtheorem{Obs}[The]{Observation}
\newtheorem{SConj}[The]{Standard Conjecture}
\newtheorem{Titre}[The]{\!\!\!\! }
\newtheorem{Conj}[The]{Conjecture}
\newtheorem{Question}[The]{Question}
\newtheorem{Prob}[The]{Problem}
\newtheorem{Def}[The]{Definition}
\newtheorem{Not}[The]{Notation}
\newtheorem{Claim}[The]{Claim}
\newtheorem{Conc}[The]{Conclusion}
\newtheorem{Ex}[The]{Example}
\newtheorem{Fact}[The]{Fact}
\newcommand{\C}{\mathbb{C}}
\newcommand{\R}{\mathbb{R}}
\newcommand{\N}{\mathbb{N}}
\newcommand{\Z}{\mathbb{Z}}
\newcommand{\Q}{\mathbb{Q}}
\newcommand{\Proj}{\mathbb{P}}

\begin{center}

{\Large\bf Deformation Openness and Closedness of Various Classes of Compact Complex Manifolds; Examples}

\end{center}

\begin{center}

 {\large Dan Popovici}

\end{center}

\vspace{1ex}

\noindent {\small {\bf Abstract.} We review the relations between compact complex manifolds carrying various types of Hermitian metrics (K\"ahler, balanced or {\it strongly Gauduchon}) and those satisfying the $\partial\bar\partial$-lemma or the degeneration at $E_1$ of the Fr\"olicher spectral sequence, as well as the behaviour of these properties under holomorphic deformations. The emphasis will be placed on the notion of {\it strongly Gauduchon} (sG) manifolds that we introduced recently in the study of deformation limits of projective and Moishezon manifolds. Various examples of sG and non-sG manifolds are exhibited while a range of constructions already known in the literature are reviewed and reinterpreted from this new standpoint.}

\vspace{3ex}

\section{Introduction}\label{section:introd}

 A Hermitian metric on a given compact complex manifold $X$ ($\mbox{dim}_{\C}X:=n$) will be identified with the associated positive-definite $C^{\infty}$ $(1, \, 1)$-form $\omega>0$ on $X$. The metric $\omega$ is said to be  \\

\noindent $\bullet$ {\it K\"ahler} if $d\omega=0$; \\

\noindent $\bullet$ {\it balanced}\footnote{terminology of [Mic82], unrelated to Donaldson's balanced metrics used in the context of cscK metrics} (or {\it semi-K\"ahler}\footnote{terminology of [Gau77b]} or {\it co-K\"ahler}) if $d\omega^{n-1}=0$; \\

\noindent $\bullet$ {\it strongly Gauduchon}\footnote{notion introduced in [Pop09]} (or {\it sG}) if $\partial\omega^{n-1}$ is $\bar\partial$-exact; \\

\noindent $\bullet$ {\it Gauduchon} if $\partial\omega^{n-1}$ is $\bar\partial$-closed or, equivalently, if $\partial\bar\partial\omega^{n-1}=0$. \\

 Although Hermitian metrics always exist, the stronger K\"ahler, balanced and sG metrics need not exist on an arbitrary $X$. A compact complex manifold carrying one of the first three types of metrics described above is said to be a compact {\it K\"ahler}, {\it balanced}, respectively {\it strongly Gauduchon} (or {\it sG}) manifold. By contrast, Gauduchon metrics exist on any compact complex manifold (cf. [Gau77a]). Actually Gauduchon's main result in [Gau77a] asserts far more\!\!: there exists a Gauduchon metric (unique up to normalisation) in the conformal class of {\it any} Hermitian metric on $X$. At the level of Hermitian metrics $\omega$ on $X$, the following obvious implications hold\!\!: \\

 $\omega$ is K\"ahler $\Longrightarrow$ $\omega$ is balanced $\Longrightarrow$ $\omega$ is sG $\Longrightarrow$ $\omega$ is Gauduchon, \\

\noindent while in the case of complex surfaces $X$ (i.e. $n=2$), it is obvious that $\omega$ is a K\"ahler metric if and only if $\omega$ is a balanced metric. An sG metric $\omega$ need not be K\"ahler even when $n=2$, but compact complex surfaces $X$ carrying sG metrics also carry K\"ahler metrics (see [Pop09, section 3]). Thus at the level of compact complex manifolds we have the following equivalences: \\

$X$ is a K\"ahler surface $\Longleftrightarrow$ $X$ is a balanced surface $\Longleftrightarrow$ $X$ is an sG surface. \\

 However, if compact complex surfaces are replaced by compact complex manifolds of complex dimension $n\geq 3$, the above equivalences break down, only strict implications from left to right hold. \\

 Recall that the so-called {\it K\"ahler currents} provide a singular (thus weaker) substitute for K\"ahler metrics. A $d$-closed positive $(1, \, 1)$-current $T$ on $X$ is said to be a {\it K\"ahler current}\footnote{terminology introduced in [JS93]} if it satisfies the strong positivity condition

$$T\geq\varepsilon\omega \hspace{3ex} \mbox{on}\,\, X,$$

\noindent for some constant $\varepsilon>0$ and some Hermitian metric $\omega>0$. K\"ahler currents need not exist on an arbitrary $X$, but they may exist when K\"ahler metrics do not. Recall that $X$ is said to be a Fujiki {\it class} ${\cal C}$ manifold if $X$ is bimeromorphic to a compact K\"ahler manifold, or equivalently, if there exists a proper holomorphic bimeromorphic map (i.e. a {\it modification}) 

$$\mu:\widetilde{X} \rightarrow X$$

\noindent from a compact K\"ahler manifold $\widetilde{X}$. Fujiki introduced {\it class} ${\cal C}$ manifolds $X$ as meromorphic images of compact K\"ahler manifolds in [Fuj78], while Varouchas gave them the above nice characterisation in [Var86]. It is a result of Demailly and Paun that {\it class} ${\cal C}$ manifolds are characterised by the existence of a K\"ahler current.

\begin{The}(Demailly-Paun [DP04])\label{The:Dem-Pau} A compact complex manifold $X$ is of class ${\cal C}$ if and only if there exists a K\"ahler current $T$ on $X$.

\end{The}

 Recall that a {\it Moishezon manifold} is, by definition, a compact complex manifold that is bimeromorphic to a projective manifold. Equivalently, $X$ is Moishezon if and only if there exists a modification 

$$\mu:\widetilde{X} \rightarrow X$$

\noindent from a projective manifold $\widetilde{X}$. Thus Moishezon manifolds are to projective manifolds what {\it class} ${\cal C}$ manifolds are to compact K\"ahler manifolds. 

 The special case of integral cohomology classes is relevant in characterisations of some of the above classes of manifolds. Recall that the De Rham cohomology $2$-class $\{\omega\}\in H^2_{DR}(X, \, \R)$ (resp. $\{T\}\in H^2_{DR}(X, \, \R)$) defined by a $C^{\infty}$ $d$-closed real $(1, \, 1)$-form $\omega$ (resp. by a $d$-closed real $(1, \, 1)$-current $T$) is said to be {\it integral} if it is the first Chern class of a holomorphic line bundle $L\rightarrow X$ or, equivalently, if $\omega$ (resp. $T$) is the curvature form (resp. curvature current) $\frac{i}{\pi}\Theta_h(L)$ of a holomorphic line bundle $(L, \, h)\rightarrow X$ endowed with a $C^{\infty}$ (resp. singular) Hermitian fibre metric $h$.

 There are neat characterisations of projective and Moishezon manifolds mirroring the general case of arbitrary (i.e. possibly transcendental) classes that occur on K\"ahler and {\it class} ${\cal C}$ manifolds.

\begin{The}(Kodaira's Embedding Theorem)\label{The:Kod-embedding} A compact complex manifold $X$ is projective if and only if there exists a K\"ahler metric $\omega$ on $X$ whose De Rham cohomology class $\{\omega\}\in H^2_{DR}(X, \, \R)$ is integral.

\end{The}

 Thus projective manifolds are {\it integral class} special cases of compact K\"ahler manifolds. Likewise, Moishezon manifolds are {\it integral class} special cases of {\it class} ${\cal C}$ manifolds as the following characterisation shows.

\begin{The}(Ji-Shiffman [JS93])\label{The:Ji-Shiff} A compact complex manifold $X$ is Moishezon if and only if there exists a K\"ahler current $T$ on $X$ whose De Rham cohomology class $\{T\}\in H^2_{DR}(X, \, \R)$ is integral. 

\end{The}

 Now recall the following standard piece of terminology.

\begin{Def}\label{Def:dd-bar-lemma} A compact complex manifold $X$ is said to satisfy the $\partial\bar\partial$-{\bf lemma} if for any $d$-closed {\it pure-type} form $u$ on $X$, the following exactness properties are equivalent: \\

$u$ is $d$-exact $\Longleftrightarrow$ $u$ is $\partial$-exact $\Longleftrightarrow$ $u$ is $\bar\partial$-exact $\Longleftrightarrow$ $u$ is $\partial\bar\partial$-exact

\end{Def}

\vspace{2ex}

 If the {\it pure-type} assumption on $u$ is dropped, then $u$ must be assumed to be both $d$-closed and $d^c$-closed (or, equivalently, both $\partial$-closed and $\bar\partial$-closed) before the above exactness properties are required to be equivalent (cf. [DGMS75]). For a {\it pure-type} form $u$, the sole $d$-closedness is equivalent to $u$ being both $\partial$-closed and $\bar\partial$-closed. However, if $u$ is not of pure-type, $du=0$ does not imply $\partial u=0$ and $\bar\partial u=0$.

 It is a standard fact in Hodge theory that any compact K\"ahler manifold satisfies the $\partial\bar\partial$-lemma. Moreover, if $\mu:\widetilde{X}\rightarrow X$ is a modification between compact complex manifolds and if the $\partial\bar\partial$-lemma holds for $\widetilde{X}$, then the $\partial\bar\partial$-lemma also holds for $X$ (see e.g. [DGMS75, Theorem 5.22.]). In particular, {\it class} ${\cal C}$ manifolds (hence also Moishezon manifolds) satisfy the $\partial\bar\partial$-lemma.

\vspace{2ex}

 An interesting result of Alessandrini and Bassanelli (see [AB91b], [AB93], [AB95]) asserts that every {\it class} ${\cal C}$ manifold is {\it balanced} (i.e. carries a {\it balanced metric}). They actually managed rather more in proving that {\it balanced manifolds} are stable under modifications.

\begin{The}\label{The:A-B-modbal}(Alessandrini-Bassanelli [AB95]) Let $\mu:\widetilde{X}\rightarrow X$ be a modification of compact complex manifolds. Then $X$ is balanced if and only if $\widetilde{X}$ is balanced.

\end{The}

 Now if $X$ is a {\it class} ${\cal C}$ manifold, then by [Var86] there exists a modification $\mu:\widetilde{X}\rightarrow X$ where $\widetilde{X}$ is a compact K\"ahler manifold. Then $\widetilde{X}$ is also balanced and, by Theorem \ref{The:A-B-modbal}, $X$ must be balanced as well. 

\vspace{2ex}

 It is worth mentioning that the sG condition enjoys the same modification stability property.

\begin{The}\label{The:mod-stab-sG}(Theorem 1.3. in [Pop10b]) Let $\mu:\widetilde{X}\rightarrow X$ be a modification of compact complex manifolds. Then $X$ is {\bf strongly Gauduchon} if and only if $\widetilde{X}$ is {\bf strongly Gauduchon}.

\end{The}

 The $\partial\bar\partial$-lemma property and the balanced property are unrelated (see e.g. examples below), but they are both implied by the {\it class} ${\cal C}$ property and, in turn, they both imply the sG property.

 Let us recall the argument given in [Pop09] proving that on every compact complex manifold $X$ which satisfies the $\partial\bar\partial$-lemma, the notions of Gauduchon and sG metrics are equivalent, hence every such $X$ is an sG manifold. If $\omega$ is a Gauduchon metric on $X$, then the pure-type $(n, \, n-1)$-form $\partial\omega^{n-1}$ is $\bar\partial$-closed by definition. Since $\partial\omega^{n-1}$ is also $\partial$-closed, it must be $d$-closed. Thus, if the $\partial\bar\partial$-lemma holds on $X$, the $d$-closed and $\partial$-exact pure-type form $\partial\omega^{n-1}$ must also be $\bar\partial$-exact. Hence $\omega$ is an sG metric on $X$.

\vspace{2ex}

 Recall that on any compact complex manifold $X$, we always have

\begin{equation}\label{eqn:b-h-eq}b_k(X)\leq\sum\limits_{p+q=k}h^{p, \, q}(X),  \hspace{3ex} k=0, 1, \dots , n=\mbox{dim}_{\C}X,\end{equation}

\noindent where $b_k(X):=\mbox{dim}_{\C}H^k_{DR}(X, \, \C)$ is the $k^{th}$ Betti number of $X$ and $h^{p, \, q}(X):=\mbox{dim}_{\C}H^{p, \, q}(X, \, \C)$ is the Hodge number of type $(p, \, q)$ of $X$. The Dolbeault cohomology group associated with $(p, \, q)$-forms on $X$ is denoted, as usual, by $H^{p, \, q}(X, \, \C)$. Recall that the Fr\"olicher spectral sequence of $X$ degenerates at $E_1$ (a property that will be denoted by $E_1(X)=E_{\infty}(X)$) if and only if equality is achieved in all of the above inequalities (\ref{eqn:b-h-eq}). Indeed, we always have $H^k_{DR}(X, \, \C)=\bigoplus\limits_{p+q=k}E^{p, \, q}_{\infty}(X)$ for all $k$ and $E_1^{p, \, q}(X)= H^{p, \, q}(X, \, \C)$ for all $p, q$. Fr\"olicher degeneration at $E_1$ means that $E_1^{p, \, q}(X)= E^{p, \, q}_{\infty}(X)$ for all $p, q$, hence it is equivalent to the existence of a Hodge decomposition on $X$, possibly {\it without} Hodge symmetry, in the following form:

$$H^k_{DR}(X, \, \C)=\bigoplus\limits_{p+q=k}H^{p, \, q}(X, \, \C),  \hspace{3ex} k=0, 1, \dots , n=\mbox{dim}_{\C}X.$$

On the other hand, recall that the $\partial\bar\partial$-lemma implies that such a Hodge decomposition, possibly {\it without} Hodge symmetry, holds on $X$. Thus the Fr\"olicher spectral sequence of any compact complex manifold $X$ that satisfies the $\partial\bar\partial$-lemma degenerates at $E_1$.

 While the sG property of $X$ and the degeneration at $E_1$ of the spectral sequence of $X$ are unrelated (see Theorem \ref{The:sG-specseq-unrel} below), they are both implied by the $\partial\bar\partial$-lemma. It is worth recalling that every compact complex {\it surface} $X$ satisfies the property $E_1(X)=E_{\infty}(X)$, while K\"ahler {\it surfaces} are the only compact complex {\it surfaces} that are sG. Thus in the case of complex {\it surfaces}, the sG property is very strong (i.e. equivalent to the K\"ahler property), while the $E_1(X)=E_{\infty}(X)$ property is trivially satisfied by all these surfaces.

\vspace{3ex}

 The relations among these properties of a compact complex manifold $X$ are summed up in the following diagram (skew arrows indicate implications)\!:

\vspace{10ex}

\noindent$\begin{array}{lllllllll}(\star) & & & & & & & & \\
          &  & X \, \mbox{Moishezon} &  & & & X \, \mbox{balanced} &  & \\
          & \rotatebox{45}{$\implies$} &  & \rotatebox{-45}{$\implies$} &  & \rotatebox{45}{$\implies$}  &   & \rotatebox{-45}{$\implies$} & \\
       X  \, \mbox{projective} &  &  & & X \, \mbox{{\it class}}\,\, {\cal C} & & & &  X  \, \mbox{sG} \\
          & \rotatebox{-45}{$\implies$} &  & \rotatebox{45}{$\implies$} &  & \rotatebox{-45}{$\implies$}   &  &  \rotatebox{45}{$\implies$} & \\
          &  & X  \, \mbox{K\"ahler} & & & & \begin{array}{l}X \,\,\mbox{satisfies}\\\partial\bar\partial-\mbox{lemma}\end{array}   &  & \\
  & & & & & & & \rotatebox{-45}{$\implies$} & \\
  & & & & & & & & E_1(X)=E_{\infty}(X)
\end{array}$

\vspace{10ex}

 One of the purposes of this work is to explain how the new notion of {\it strongly Gauduchon manifold} fits into the context of the well-known earlier other notions featuring in diagram $(\star)$ above. On the one hand, we shall provide examples of compact complex manifolds that are not sG (cf. section \ref{section:not-sG}) in order to stress that the sG condition is not automatically satisfied. 

\begin{The}\label{The:not-sG} The Calabi-Eckmann manifolds [CE53], the Hopf manifolds [Hop48] and Tsuji's manifolds [Tsu84] are {\bf not strongly Gauduchon}.

\end{The}

 On the other hand, we shall provide examples of compact complex manifolds that 
drive home the peculiarity of the sG condition distinguishing it from the two stronger notions that immediately precede it in diagram $(\star)$ (cf. also Theorem \ref{The:sG-examp} in section \ref{section:sG-examp}).

\begin{The}\label{The:sG-examp} There exist compact complex manifolds that are {\bf strongly Gauduchon} but are {\bf not balanced} and on which the {\bf $\partial\bar\partial$-lemma does not hold}.

\end{The}

 We also point out that the last two notions in diagram $(\star)$ are unrelated.

\begin{The}\label{The:sG-specseq-unrel} 

$(a)$\, There exist {\bf strongly Gauduchon} compact complex manifolds whose {\bf Fr\"olicher spectral sequence does not degenerate at $E_1$} (e.g. the Iwasawa manifold, which is even balanced). 

$(b)$\, There exist compact complex manifolds whose {\bf Fr\"olicher spectral sequence degenerates at $E_1$} but which are {\bf not strongly Gauduchon} (e.g. any compact non-K\"ahler complex surface).

\end{The}

 This observation seems to indicate the existence of two disjoint realms in diagram $(\star)$: the {\it metric} notions (i.e. {\it balanced} and {\it sG}) branching off upwards from the {\it class} ${\cal C}$ notion; the {\it topological} notions (i.e. $\partial\bar\partial$-{\it lemma} and $E_1=E_{\infty}$) branching off downwards. The overall idea underlying the proof of our main result in [Pop09] was based on switching from the latter to the former realm: after observing that the {\it topological} condition of non-jumping at $t=0$ of the Hodge number $h^{0, \, 1}(t)$ is sufficient for our purposes but next to impossible to guarantee on an {\it a priori} basis, we replaced it with the {\it metric} sG condition on the central fibre that we managed to guarantee.    

 Another purpose of this work is to review known results and open questions about the behaviour of these properties under holomorphic deformations with a special emphasis on the new notion of {\it strongly Gauduchon} (or sG) manifolds. We adopt the Kodaira-Spencer terminology (cf. e.g. [Kod86]): \\

{\it A complex analytic (or holomorphic) family of compact complex manifods is a proper holomorphic submersion $\pi:{\cal X}\to\Delta$ from an arbitrary complex analytic manifold ${\cal X}$ to some open ball $\Delta\subset\C^m$ about the origin $0\in\C^m$.} \\

  Thus the fibres $X_t:=\pi^{-1}(t)$, $t\in\Delta$, are all (smooth) compact complex manifolds of the same dimension varying holomorphically with the parameter $t\in\Delta$. It is well-known that any such family is differentiably trivial, i.e. there exists a $C^{\infty}$ manifold $X$ independent of $t\in\Delta$ such that the fibre $X_t$ is $C^{\infty}$-diffeomorphic to $X$ for all $t\in\Delta$. Only the complex structure $J_t$ of $X_t$ varies with $t\in\Delta$. Thus the holomorphic family $(X_t)_{t\in\Delta}$ can be identified with a fixed $C^{\infty}$ manifold $X$ endowed with a holomorphic family of complex structures $(J_t)_{t\in\Delta}$ (cf. [Kod86]). 

 It will be sufficient to restrict attention to the case where the base $\Delta$ is an open disc about the origin in $\C$, i.e. $m=1$. We shall be concerned with stability properties of the notions featuring in diagram $(\star)$ when they appear in holomorphic families as above. Two points of view will be adopted.

\begin{Def} $(i)$\, A given property $(P)$ of a compact complex manifold is said to be {\bf open} under holomorphic deformations if for every holomorphic family of compact complex manifolds $(X_t)_{t\in\Delta}$ and for every $t_0\in\Delta$, the following implication holds: \\

$X_{t_0}$ has property $(P)$ $\implies$  $X_t$ has property $(P)$ for all $t\in\Delta$ sufficiently close to $t_0.$ \\

\noindent $(ii)$\, A given property $(P)$ of a compact complex manifold is said to be {\bf closed} under holomorphic deformations if for every holomorphic family of compact complex manifolds $(X_t)_{t\in\Delta}$ and for every $t_0\in\Delta$, the following implication holds: \\

$X_t$ has property $(P)$ for all $t\in\Delta\setminus\{t_0\}$ $\implies$  $X_{t_0}$ has property $(P).$

\end{Def}

\vspace{2ex}

 The interest in deformation stability questions was sparked by the following celebrated result of Kodaira and Spencer.

\begin{The}(Kodaira-Spencer [KS60])\label{The:KS-open} The {\bf K\"ahler property} of compact complex manifolds is {\bf open} under holomorphic deformations.

\end{The}

 However, when Hironaka constructed an example in 1962, the following fact came as a bit of a surprise.

\begin{The}(Hironaka [Hir62])\label{The:Hir-notclosed} The {\bf K\"ahler property} of compact complex manifolds of complex dimension $\geq 3$ is {\bf not closed} under holomorphic deformations.

\end{The}

 Our main result of [Pop09] states that the degeneration of the projective property of compact complex manifolds in the deformation limit is relatively mild, i.e. projective manifolds degenerate to Moishezon manifolds. In view of Hironoka's example that proved Theorem \ref{The:Hir-notclosed}, the following result is optimal.

\begin{The}(Theorem 1.1. in [Pop09])\label{The:proj-def} Let $\pi:{\cal X}\to\Delta$ be a complex analytic family of compact complex manifolds. If the fibre $X_t$ is projective for every $t\in\Delta^{\star}$, then $X_0$ is Moishezon.

\end{The}

 Hironaka's example does not cover the case of holomorphic families of compact complex surfaces which turned out to behave very differently from manifolds of higher dimensions. Indeed, Kodaira's classification of surfaces, Miyaoka's result [Miy74] asserting that an elliptic surface is K\"ahler if and only if its first Betti number is even and Siu's result [Siu83] asserting that every K3 surface is K\"ahler showed that the K\"ahlerness of compact complex surfaces is a topological property\!\!:

\begin{The}\label{The:surfaces_b1even}(Kodaira, Miyaoka, Siu) A compact complex {\bf surface} $X$ is K\"ahler if and only if its first Betti number $b_1(X)$ is even.

\end{The}

 Direct proofs of this theorem which do not invoke Kodaira's classification of compact complex surfaces were subsequently given by Buchdahl [Buc99] and Lamari [Lam99] independently. The reader will find further details on the history of the result in these references.

 Since all the fibres $X_t$ in a holomorphic family of compact complex manifolds $(X_t)_{t\in\Delta}$ are differentiably diffeomorphic, they have the same Betti numbers. In particular, we have

\begin{Cor}\label{Cor:K-closed-surf} The {\bf K\"ahler property} of compact complex {\bf surfaces} is both {\bf open} and {\bf closed} under holomorphic deformations. 
 
 In particular, if in a holomorphic family of compact complex {\bf surfaces} some fibre is K\"ahler, then all the fibres are K\"ahler.

\end{Cor}

 Let us also recall that the deformation behaviour of the {\it class} ${\cal C}$ property is the opposite of that of the K\"ahler property.

\begin{The}(Campana [Cam91a], Lebrun-Poon [LP92])\label{The:classC-nop} The {\bf class} ${\cal C}$ {\bf property} of compact complex manifolds is {\bf not open} under holomorphic deformations.

\end{The}  

 The examples proving this statement use families of twistor spaces and will be alluded to in section \ref{section:ES}. Since it is known by a result of Campana [Cam91b] that the Moishezon and {\it class} ${\cal C}$ properties are equivalent for twistor spaces, the Moishezon property is also seen not to be open under deformations. This latter fact is hardly surprising since a property associated with integral classes is not naturally expected to be deformation open. However, the following long-standing conjecture is still open.

\begin{SConj}\label{Conj:classC-closed} The {\bf class} ${\cal C}$ {\bf property} of compact complex manifolds is {\bf closed} under holomorphic deformations.

\end{SConj}

 Our main result in [Pop10a] amounts to a confirmation of this conjecture in the {\it integral class} case since Moishezon manifolds can be seen as integral class versions of {\it class} ${\cal C}$ manifolds (thanks to Theorems \ref{The:Dem-Pau} and \ref{The:Ji-Shiff}).

\begin{The}(Theorem 1.1. in [Pop10a])\label{The:Moi-closed} The {\bf Moishezon property} of compact complex manifolds is {\bf closed} under holomorphic deformations.

\end{The}

 We have shown in [Pop10a] (see also Theorem \ref{The:sG-openness} below) that the {\bf sG} property is deformation {\bf open}. We hope that the following also holds.

\begin{Conj}\footnote{This conjecture was suggested to the author by Jean-Pierre Demailly.}\label{Conj:sG-closedness} The {\bf sG} property of compact complex manifolds is {\bf closed} under holomorphic deformations.

\end{Conj}

 The proof of this fact is still elusive, but if this turns out to be the case, the sG property would be the only known property of compact complex manifolds to be stable under all known operations (i.e. both open and closed under deformations, as well as stable under modifications as shown in [Pop10b, cf. Theorem \ref{The:mod-stab-sG} above]). In particular, it would suffice for just one fibre $X_{t_0}$ to be sG in order to guarantee that all the fibres $X_t$ are sG. We have proved the conclusion of Conjecture \ref{Conj:sG-closedness} under a stronger assumption on the generic fibres and this has played a major part in our proofs of the main results of [Pop09] and [Pop10a].

\begin{Prop}(Proposition 4.1. in [Pop09])\label{Prop:dd-bar-sG} Let $\pi:{\cal X}\to\Delta$ be a complex analytic family of compact complex manifolds. If the $\partial\bar\partial$-lemma holds on $X_t$ for every $t\in\Delta^{\star}$, then $X_0$ is a {\bf strongly Gauduchon} manifold.

\end{Prop}

 When it comes to balanced manifolds, Alessandrini and Bassanelli showed in [AB90] (see Theorem \ref{The:bal-notopen} below) that the {\bf balanced} property of compact complex manifolds is {\bf not} deformation {\bf open}. This difference in deformation behaviour between sG and balanced manifolds will be exploited in section \ref{section:sG-examp}. However, we hope that the balanced analogue of Conjecture \ref{Conj:sG-closedness} holds.

\begin{Conj}\label{Conj:bal-closedness} The {\bf balanced} property of compact complex manifolds is {\bf closed} under holomorphic deformations.

\end{Conj}

 If this turns out to be the case, Conjecture \ref{Conj:bal-closedness} might be used to tackle the standard Conjecture \ref{Conj:classC-closed}. Indeed, by the Alessandrini-Bassanelli theorem \ref{The:A-B-modbal}, proving that the limit fibre $X_0$ is balanced when $X_t$ has been supposed to be {\it class} ${\cal C}$ (hence also balanced) for all $t\neq 0$, is necessary to proving the stronger {\it class} ${\cal C}$ property of $X_0$.

\vspace{2ex}

 As for the property of Fr\"olicher degeneration at $E_1$, we have

\begin{The}\label{The:E1=Einf_op-cl} $(a)$\, (Kodaira-Spencer [KS60]) For compact complex manifolds, the property of the {\bf Fr\"olicher spectral sequence degenerating at $E_1$} is {\bf open} under holomorphic deformations. 

$(b)$\, (Eastwood-Singer [ES93, Theorem 5.4.]) For compact complex manifolds, the property of the {\bf Fr\"olicher spectral sequence degenerating at $E_1$} is {\bf not closed} under holomorphic deformations.

\end{The}

 Part $(a)$ is by now a classical statement that follows immediately from inequality (\ref{eqn:b-h-eq}) satisfied by every fibre $X_t$ (in which equality is equivalent to $E_1(X_t)=E_{\infty}(X_t)$), from the Betti numbers $b_k(X_t)$ of the fibres being independent of $t$ (thanks to the $C^{\infty}$ triviality of the family) and from the upper-semicontinuity of every Hodge number $h^{p, \, q}(X_t)$ w.r.t. $t\in\Delta$ (which is another classical result of Kodaira and Spencer). An overview of the Eastwood-Singer proof of part $(b)$ will be given in section \ref{section:ES}.

 We are at a loss to know anything about the deformation properties of compact complex manifolds satisfying the $\partial\bar\partial$-lemma. It might well be the case that this property is neither open nor closed under deformations although the evidence is very tenuous. We will outline a possible approach to this question via twistor spaces in the last section \ref{section:ES}.

\vspace{2ex}

\noindent {\bf Acknowledgments.} The author is extremely grateful to Professor Akira Fujiki for kindly inviting him to Osaka University, for patiently explaining to him exciting titbits about various notions in mathematics and for indicating a host of bibliographical references that have broadened his understanding of a wider picture. Thanks are also due to Professor Hajime Tsuji over a similar invitation to Tokyo and for kindly pointing out the reference [Tsu84].

\section{Examples of non-sG compact complex manifolds}\label{section:not-sG}

 In this section we prove Theorem \ref{The:not-sG} by exhibiting three well-known classes of compact complex manifolds that are not sG\!\!: the Calabi-Eckmann manifolds, the Hopf manifolds and Tsuji's manifolds constructed in [Tsu84]. The underlying space of all these manifolds is a product $X:=S^{2p+1}\times S^{2q+1}$ of two real odd-dimensional spheres, so they all share the property $H^2_{DR}(X, \, \R)=0$ for the second De Rham cohomology group. This implies that any $d$-closed positive current $T$ of type $(1, \, 1)$ on $X$, should it exist, must be $d$-exact since the associated De Rham cohomology $2$-class $\{T\}\in H^2_{DR}(X, \, \R)$ must vanish. The existence of a non-trivial $(1, \, 1)$-current $T$ on $X$ that is both positive and $d$-exact amounts to $X$ being non-sG as the following intrinsic characterisation of sG manifolds obtained in [Pop09] shows.

\begin{Prop}\label{Prop:current-sG-char}(Proposition 3.3. in [Pop09]) Let $X$ be a compact complex manifold. Then $X$ carries a strongly Gauduchon metric if and only if there exists no non-zero current $T$ of type $(1, \, 1)$ such that $T\geq 0$ and $T$ is $d$-exact on $X$.

\end{Prop}

 We shall briefly review the three classes of compact complex manifolds mentioned above and notice that every such manifold $X$ possesses complex hypersurfaces $Y\subset X$. Thus, since $H^2_{DR}(X, \, \R)=0$, the current of integration on any of these complex hypersurfaces $Y$ is a current as in Proposition \ref{Prop:current-sG-char}, ruling out the possibility that any manifold $X$ in one of these classes be sG.

\vspace{2ex}

\noindent $(a)$\, {\bf Calabi-Eckmann manifolds.} For all $p, q\in\N$, Calabi and Eckmann [CE53] constructed a complex structure on the Cartesian product $S^{2p+1}\times S^{2q+1}$ of odd-dimensional spheres. The case $p=q=0$ being equivalent to a closed Riemann surface of genus $1$ and periods $1, \tau$, they assume $p>0$. In the case $q=0$, the Calabi-Eckmann complex structure on $S^{2p+1}\times S^1$, although constructed by a different method, coincides with the complex structure constructed earlier by Hopf in [Hop48] starting from the universal covering space of $S^{2p+1}\times S^1$ equipped with the complex structure of $\C^{p+1}\setminus\{0\}$. The simply connected manifolds $S^{2p+1}\times S^{2q+1}$ $(p, q>0)$ are given in [CE53] complex structures making them into compact, simply connected, non-K\"ahler complex manifolds $M^{p, \, q}$ of complex dimension $p+q+1$ enjoying, among others, the following properties (for all $p, q$, including $q=0$)\!\!: \\

$(i)$\, there exists a complex analytic fibring $\sigma : M^{p, \, q}\to\Proj^p\times\Proj^q$ over the product of complex projective spaces $\Proj^p$ and $\Proj^q$ whose fibres are tori of real dimension $2$ (or algebraic curves of genus $1$) (cf. [CE56, Theorem II]); \\

$(ii)$\, every compact complex subvariety of $M^{p, \, q}$ is the set of all points that are mapped by $\sigma$ onto an algebraic subvariety of $\Proj^p\times\Proj^q$; it is therefore also fibred by tori (cf. [CE56, Theorem IV]). \\

 It is clear that the inverse image under $\sigma$ of any complex hypersurface of $\Proj^p\times\Proj^q$ defines a complex hypersurface of the Calabi-Eckmann manifold $M^{p, \, q}$. Thus no Calabi-Eckmann manifold $M^{p, \, q}$ ($p>0$) can be an sG manifold.\footnote{This same argument was invoked in [Mic82, p.263] to show that Calabi-Eckmann manifolds are not balanced.}

\vspace{2ex}

\noindent $(b)$\, {\bf Hopf manifolds.} As mentioned above (and proved in $\S.3$ of [CE56]), the Hopf manifolds $S^{2p+1}\times S^1$ ($p>0$) endowed with the complex structure constructed in [Hop48] can be seen in retrospect as special cases for $q=0$ of Calabi-Eckmann manifolds. Thus they contain complex hypersurfaces and are not sG manifolds by the above arguments.

\vspace{2ex}

\noindent $(c)$\, {\bf Tsuji's manifolds.} Generalising the Calabi-Eckmann complex structures, Tsuji constructed in [Tsu84] complex structures on $S^3\times S^3$ in the following way. Starting from an arbitrary $(\alpha_1, \alpha_2, \alpha_3)\in\C^3$ satisfying

$$0<|\alpha_1|\leq |\alpha_2|<1 \hspace{2ex} \mbox{and} \hspace{2ex} 0<|\alpha_3|<1,$$

\noindent the author of [Tsu84] considers the primary Hopf manifold of complex dimension $3$ 

$$H(\alpha):=\C^3\setminus\{0\}/\langle h\rangle,$$

\noindent where the automorphism $h:\C^3\to\C^3$ is defined by $h(z_1, z_2, z_3):=(\alpha_1\, z_1, \, \alpha_2\, z_2, \, \alpha_3\, z_3)$ for all $(z_1, z_2, z_3)\in\C^3$ and $\langle h\rangle\subset\mbox{Aut}(\C^3)$ denotes the automorphism group generated by $h$. He then goes on to consider

$$C:=\{[z_1, z_2, z_3]\in H(\alpha)\,\, ; \,\, z_1=z_2=0\}\subset H(\alpha),$$

\noindent an elliptic curve contained in $H(\alpha)$ and

$$S_0:=\{[z_1, z_2, z_3]\in H(\alpha)\,\, ; \,\, z_3=0\}\subset H(\alpha),$$

\noindent a primary Hopf surface which is a complex hypersurface of $H(\alpha)$. For every 

$$A=\begin{pmatrix}a & b\\
                      c & d\end{pmatrix}\in SL(2, \, \Z) \hspace{2ex} \mbox{and} \hspace{2ex} m=(m_1, \, m_2)\in\Z^2, m_1, m_2\gg 1,$$

\noindent the author shows the existence of $\beta=(\beta_1, \beta_2, \beta_3)\in\C^3$ defining biholomorphisms

$$L^{\star}(\beta)\stackrel{\Phi^{\pm}}{\simeq} L^{\star}(\alpha),$$

\noindent where $L^{\star}(\alpha)$ and $L^{\star}(\beta)$ are obtained from $L(\alpha)$ and $L(\beta)$ by removing the respective zero section, while $L(\alpha)$ and $L(\beta)$ are holomorphic line bundles over the respective primary Hopf surfaces

$$S_{\alpha_1, \, \alpha_2, \, 0}:=\C^2\setminus\{0\}/\langle g_{\alpha}\rangle  \hspace{2ex} \mbox{and} \hspace{2ex} S_{\beta_1, \, \beta_2, \, 0}:=\C^2\setminus\{0\}/\langle g_{\beta}\rangle$$

\noindent associated with automorphisms of $\C^2$

$$g_{\alpha}(z_1, \, z_2):=(\alpha_1\, z_1, \, \alpha_2\, z_2)  \hspace{2ex} \mbox{and} \hspace{2ex} g_{\beta}(z_1, \, z_2):=(\beta_1\, z_1, \, \beta_2\, z_2)$$

\noindent defined by 

$$L(\alpha):=\C^2\setminus\{0\}\times\C/\langle h_{\alpha}\rangle \hspace{2ex} \mbox{and} \hspace{2ex} L(\beta):=\C^2\setminus\{0\}\times\C/\langle h_{\beta}\rangle,$$

\noindent where the automorphisms $h_{\alpha}$ and $h_{\beta}$ of $\C^3$ are defined by

$$h_{\alpha}(z_1, \, z_2, \, z_3):=(\alpha_1\, z_1, \, \alpha_2\, z_2, \alpha_3\, z_3)  \hspace{2ex} \mbox{and} \hspace{2ex} h_{\beta}(z_1, \, z_2, \, z_3):=(\beta_1\, z_1, \, \beta_2\, z_2, \, \beta_3\, z_3).$$

 Considering a compactification of $L(\beta)$ as a $\Proj^1$-bundle $\Proj(\beta)\to S_{\beta_1, \, \beta_2, \, 0}$, the infinity section of $\Proj(\beta)$ is denoted $S_{\infty}$, while $U(S_{\infty})$ denotes a tubular neighbourhood of $S_{\infty}$ in $\Proj(\beta)$. The author defines compact complex manifolds 

$$M^{\pm}(\alpha, A, m)$$

\noindent by identifying 

$$L^{\star}(\beta)\subset\Proj(\beta)\setminus(\mbox{zero section})$$

\noindent with

 $$L^{\star}(\alpha)\simeq H(\alpha)\setminus (S_0\cup C)\subset H(\alpha)$$

\noindent using $\Phi^{\pm}$. These compact complex manifolds are seen to arise as

\begin{equation}\label{eqn:surgery}M^{\pm}(\alpha, A, m)= (H(\alpha)\setminus C)\cup U(S_{\infty}),\end{equation}

\noindent or equivalently, $M^{\pm}(\alpha, A, m)$ are obtained from $H(\alpha)$ by a surgery which replaces $C$ with $U(S_{\infty})$.

\begin{The}([Tsu84, Theorem 1.13]) $M^{\pm}(\alpha, A, m)$ is diffeomorphic to $S^3\times S^3$ if and only if $A$ is of the form $A=\begin{pmatrix}a & b\\
                                                                   \pm 1 & d\end{pmatrix}.$

 Consequently, if $A$ has the above shape, $M^{\pm}(\alpha, A, m)$ is diffeomorphic to an $S^3$-bundle over a lens space, hence $M^{\pm}(\alpha, A, m)$ has a complex structure.

\end{The}

 With this outline of Tsuji's construction understood, we see that the complex hypersurface $S_0\subset H(\alpha)$ satisfies $S_0\cap C = \emptyset$. Thus, in view of the description (\ref{eqn:surgery}) of $M^{\pm}(\alpha, A, m)$, we get a complex hypersurface

$$S_0\subset M^{\pm}(\alpha, A, m)$$

\noindent whose existence, along with the property $H^2_{DR}(M^{\pm}(\alpha, A, m), \, \R)=0$, shows that Tsuji's compact complex manifolds $M^{\pm}(\alpha, A, m)$ are not sG for any $\alpha\in\C^3, A\in SL(2, \, \Z), m=(m_1, \, m_2)\in\Z^2$ as above.

\section{Examples of sG manifolds}\label{section:sG-examp}

 As pointed out in the Introduction, all balanced manifolds and all compact complex manifolds on which the $\partial\bar\partial$-lemma holds provide examples of sG manifolds (cf. implication diagram $(\star)$). In this section we shall exhibit compact complex manifolds that are sG but neither are they balanced nor do they satisfy the $\partial\bar\partial$-lemma. These examples are thus meant to emphasise the difference between sG manifolds and stronger earlier types of possibly non-K\"ahler compact complex manifolds.

 The examples we shall exhibit will be provided by holomorphic families of compact complex manifolds. The starting point is the following stability property of $sG$ manifolds under small deformations.

\begin{The}([Pop10a, Conclusion 2.4.])\label{The:sG-openness} The {\bf sG property} of compact complex manifolds is {\bf open} under holomorphic deformations.

\end{The}

 We now recall the proof of this {\it small deformation} stability result. The main argument is provided by the following characterisation of sG manifolds.

\begin{Lem}([Pop09, Lemma 3.2])\label{Lem:sGposchar} Let $X$ be a compact complex manifold of complex dimension $n$. Then $X$ carries an sG metric if and only if there exists a $C^{\infty}$ $(2n-2)$-form $\Omega$ on $X$ satisfying the following three conditions: \\

\noindent $(a)$\, $\Omega=\overline{\Omega}$ (i.e. $\Omega$ is real);

\noindent $(b)$\, $d\Omega=0$;

\noindent $(c)$\, $\Omega^{n-1, \, n-1}>0$ on $X$ (i.e. the component of type $(n-1, \, n-1)$ of $\Omega$ w.r.t. the complex structure of $X$ is positive-definite).

\end{Lem}

 Note that conditions $(a)$ and $(b)$ are independent of the complex structure of $X$, while a change of complex structure changes the $(n-1, \, n-1)$-component of a given $(2n-2)$-form $\Omega$. Thus condition $(c)$ is the only one to be dependent on the complex structure of $X$.

\vspace{2ex}

\noindent {\it Proof of Lemma \ref{Lem:sGposchar}(cf. [Pop09].)} The vanishing of the $(2n-1)$-form $d\Omega$ (cf. $(b)$) amounts to the simultaneous vanishing of its components $\partial\Omega^{n-1, \, n-1} + \bar\partial\Omega^{n, \, n-2}$ (of type $(n, \, n-1)$) and $\partial\Omega^{n-2, \, n} + \bar\partial\Omega^{n-1, \, n-1}$ (of type $(n-1, \, n)$). These two components are conjugate to each other if $\Omega$ satisfies $(a)$. Thus, if $(a)$ holds, $(b)$ is equivalent to $\partial\Omega^{n-1, \, n-1} + \bar\partial\Omega^{n, \, n-2}=0.$

 Suppose there exists an sG metric $\omega$ on $X$. This means that $\omega$ is a $C^{\infty}$ positive-definite $(1, \, 1)$-form on $X$ such that the $(n, \, n-1)$-form $\partial\omega^{n-1}$ is $\bar\partial$-exact. Then the $(n-1, \, n-1)$-form $\Omega^{n-1, \, n-1}:=\omega^{n-1}$ is positive-definite on $X$ and there exists a $C^{\infty}$ $(n, \, n-2)$-form $\Omega^{n, \, n-2}$ on $X$ satisfying $\partial\Omega^{n-1, \, n-1}=-\bar\partial\Omega^{n, \, n-2}$. Considering the $(n-2, \, n)$-form $\Omega^{n-2, \, n}:=\overline{\Omega^{n, \, n-2}}$, we see that the $C^{\infty}$ $(2n-2)$-form

$$\Omega:=\Omega^{n, \, n-2} + \Omega^{n-1, \, n-1} + \Omega^{n-2, \, n}$$

\noindent satisfies conditions $(a)$, $(b)$, $(c)$.  

 Conversely, suppose there exists a $C^{\infty}$ $(2n-2)$-form $\Omega$ on $X$ satisfying conditions $(a)$, $(b)$, $(c)$. According to an observation in linear algebra due to Michelsohn [Mic82], every $C^{\infty}$ positive-definite $(n-1, \, n-1)$-form admits a unique $(n-1)^{st}$ root. Applying this to $\Omega^{n-1, \, n-1}>0$, we get a unique $C^{\infty}$ positive-definite $(1, \, 1)$-form $\omega>0$ on $X$ such that 

$$\omega^{n-1}=\Omega^{n-1, \, n-1}.$$

\noindent By condition $(b)$ satisfied by $\Omega$, we see that $\partial\omega^{n-1}$ is $\bar\partial$-exact, which means that the Hermitian metric $\omega$ of $X$ is {\it strongly Gauduchon}.   \hfill $\Box$

\vspace{3ex}

\noindent {\it Proof of Theorem \ref{The:sG-openness}.} If we are given a holomorphic family of compact complex manifolds $(X_t)_{t\in\Delta}$ with $\mbox{dim}_{\C}X_t=n$, we denote by $X$ the $C^{\infty}$ manifold underlying the fibres $X_t$ and by $J_t$ the complex structure of $X_t$ for all $t\in\Delta$. The family $(X_t)_{t\in\Delta}$ is thus equivalent to the holomorphic family of complex structures $(J_t)_{t\in\Delta}$ on $X$. If we have a $C^{\infty}$ $(2n-2)$-form $\Omega$ on $X$, its components $\Omega_t^{n-1, \, n-1}$ of type $(n-1, \, n-1)$ w.r.t. the complex structures $J_t$ vary in a $C^{\infty}$ way with $t\in\Delta$. Consequently, if $\Omega_0^{n-1, \, n-1}>0$ then $\Omega_t^{n-1, \, n-1}>0$ for $t\in\Delta$ sufficiently close to $0\in\Delta$. Thus condition $(c)$ of Lemma \ref{Lem:sGposchar} is preserved under small deformations by mere continuity. Since conditions $(a)$ and $(b)$ of Lemma \ref{Lem:sGposchar} are independent of the complex structure of $X$, it follows that any $C^{\infty}$ $(2n-2)$-form $\Omega$ on $X$ satisfying conditions $(a)$, $(b)$ and $(c)$ of Lemma \ref{Lem:sGposchar} w.r.t. $J_0$ also satisfies these conditions w.r.t. $J_t$ for all $t$ sufficiently near $0$. The proof of Theorem \ref{The:sG-openness} is complete.  \hfill $\Box$

\vspace{3ex}

 One fundamental difference between balanced and sG manifolds that we shall exploit is that, unlike sG manifolds, balanced manifolds are not stable under small deformations. This result was first observed by Alessandrini and Bassanelli [AB90] and refutes Michelsohn's claim of the contrary made in the introduction to [Mic82].

\begin{The}(Alessandrini-Bassanelli [AB90])\label{The:bal-notopen} The {\bf balanced} property of compact complex manifolds is {\bf not open} under holomorphic deformations.

\end{The}

 Alessandrini and Bassanelli use the explicit description of the Kuranishi family of the Iwasawa manifold (known to be balanced) calculated by Nakamura in [Nak75] and observe that one particular direction among the six dimensions of the base space yields the example that proves Theorem \ref{The:bal-notopen}. We now review a few basic facts about compact {\it complex parallelisable} manifolds to the class of which the Iwasawa manifold belongs, before surveying the arguments and results of [Nak75], [AB90] and [AB91a] that are necessary to the understanding of the observation of Alessandrini and Bassanelli (cf. Proposition \ref{Prop:A-B_examp}) which proves Theorem \ref{The:bal-notopen}.

\vspace{2ex}

 Let $X$ be a compact complex manifold, $\mbox{dim}_{\C}X=n$. Recall the following \\

\noindent {\bf Standard fact.} {\it If $X$ is K\"ahler, then for every $p = 0, 1, \dots , n$ and for every form $u\in C^{\infty}(X, \, \Lambda^{p, \, 0}T^{\star}X)$ such that $\bar\partial u=0$, we have $du=0$. 

 In other words, every holomorphic $p$-form is $d$-closed on a compact K\"ahler manifold.} \\

 To see this, fix any K\"ahler metric $\omega$ on $X$ and recall that the associated $d$-Laplacian $\Delta_{\omega}:=dd^{\star} + d^{\star}d$, $\partial$-Laplacian $\Delta'_{\omega}:=\partial\partial^{\star} + \partial^{\star}\partial$ and $\bar\partial$-Laplacian $\Delta''_{\omega}:=\bar\partial\bar\partial^{\star} + \bar\partial^{\star}\bar\partial$ (where the formal adjoints $d^{\star}$, $\partial^{\star}$, $\bar\partial^{\star}$ are calculated w.r.t. $\omega$) are related by

\begin{equation}\label{eqn:Klaplacians-eq}\Delta'_{\omega}=\Delta''_{\omega}=\frac{1}{2}\Delta_{\omega}.\end{equation}

\noindent This property is peculiar to K\"ahler metrics; it fails dramatically for arbitrary, non-K\"ahler metrics. Now, given any smooth $(p, \, 0)$-form $u$ on $X$, we clearly have $\bar\partial^{\star}u=0$ for trivial bidegree reasons. If $u$ is holomorphic (i.e. $\bar\partial u=0$), then $\Delta''_{\omega}u=0$ because $\ker\Delta''_{\omega}=\ker\bar\partial\cap\ker\bar\partial^{\star}$. (This last identity of kernels is valid for any Hermitian metric $\omega$ on a {\it compact} manifold.) It then follows from (\ref{eqn:Klaplacians-eq}) that $\Delta_{\omega}u=0$. Since $\ker\Delta_{\omega}=\ker d\cap\ker d^{\star}$, we see that $d u=0$.   \hfill $\Box$

\vspace{2ex}

  Let us now recall that this standard fact enables one to see that the Iwasawa manifold is not K\"ahler. The {\bf Iwasawa manifold} is, by definition, the compact complex manifold of complex dimension $3$ defined as the quotient

$$X:=G/\Gamma$$

\noindent of the simply connected, connected {\it complex} Lie group (the Heisenberg group)

$$G:=\left\{\begin{pmatrix}1 & z_1 & z_3\\   
                           0 & 1 & z_2\\
                   0 & 0 & 1\end{pmatrix}\,\, ; \,\, z_1, z_2, z_3\in\C\right\}\subset Gl_3(\C)$$

\noindent by the discrete subgroup $\Gamma\subset G$ of matrices with entries $z_1, z_2, z_3\in\Z[i]$. The complex manifold structure on $G$ is defined by the complex structure of $\C^3$ via the obvious diffeomorphism $G\simeq\C^3$, while the group structure on $G$ is defined by the multiplication of matrices

$$\begin{pmatrix}1 & z_1 & z_3\\   
                           0 & 1 & z_2\\
                   0 & 0 & 1\end{pmatrix}\, \begin{pmatrix}1 & w_1 & w_3\\   
                           0 & 1 & w_2\\
                   0 & 0 & 1\end{pmatrix} = \begin{pmatrix}1 & z_1+w_1 & z_3+w_3\\   
                           0 & 1 & z_2+w_2\\
                   0 & 0 & 1\end{pmatrix}.$$

\noindent Since the holomorphic $1$-form on $G$

$$G\ni M\mapsto M^{-1}\, dM$$

\noindent is invariant under the action of $\Gamma$, it descends to a holomorphic $1$-form on $X$. An elementary calculation shows that

$$\mbox{if} \hspace{2ex} M=\begin{pmatrix}1 & z_1 & z_3\\   
                           0 & 1 & z_2\\
                   0 & 0 & 1\end{pmatrix} \hspace{2ex} \mbox{then} \hspace{2ex}  M^{-1}\, dM = \begin{pmatrix}0 & dz_1 & dz_3-z_1\, dz_2\\   
                           0 & 0 & dz_2\\
                   0 & 0 & 0\end{pmatrix}.$$

\noindent Thus we get holomorphic $1$-forms on the Iwasawa manifold $X$ induced by the following forms on $\C^3$\!\!:

\begin{equation}\label{eqn:Iwa-forms}\varphi_1:=dz_1, \hspace{2ex} \varphi_2:=dz_2, \hspace{2ex} \varphi_3:=dz_3-z_1dz_2.\end{equation}

\noindent Denoting the induced forms by the same symbols $\varphi_1, \varphi_2, \varphi_3$, it is obvious that 

\begin{equation}\label{eqn:Iwa-1forms-rel}d\varphi_1=d\varphi_2=0 \hspace{2ex} \mbox{but}\,\, d\varphi_3=-\varphi_1\wedge\varphi_2\neq 0 \hspace{2ex} \mbox{on}\,\, X.\end{equation}

\noindent Since the holomorphic $1$-form $\varphi_3$ on $X$ is not $d$-closed, $X$ is not K\"ahler.

 Finally recall that the forms defined in (\ref{eqn:Iwa-forms}), which are linearly independent at every point of $X$, can be used to completely calculate the De Rham, Dolbeault and Bott-Chern cohomologies of the Iwasawa manifold (see e.g. [Sch07, p. 4-6] for details). For example, $\varphi_1$, $\varphi_2$, $\overline{\varphi_1}$, $\overline{\varphi_2}$ are all $d$-closed, but not $d$-exact, $1$-forms on $X$ and any two of them are not $d$-cohomologous. For instance, if we had $\varphi_1=df=\partial f + \bar\partial f$ for a smooth function $f$ on $X$, then $\bar\partial f=0$ since $\varphi_1$ is of type $(1, \, 0)$. Thus $f$ would be holomorphic, hence constant, on the compact $X$, which is impossible since $\varphi_1$ is not zero. It is readily seen that $H^1_{DR}(X, \, \C)$ is the $\C$-vector space generated as follows\!\!:

\begin{equation}\label{eqn:H1-Iwa}H^1_{DR}(X, \, \C) = \langle\{\varphi_1\}, \, \{\varphi_2\}, \, \{\overline{\varphi_1}\}, \, \{\overline{\varphi_2}\}\rangle,  \hspace{2ex} \mbox{hence}\hspace{1ex} b_1(X)=4,\end{equation}

\noindent where $\{\,\,\,\,\}$ denotes a De Rham cohomology class. Using (\ref{eqn:Iwa-1forms-rel}), we can easily see that $\varphi_1\wedge\varphi_3$ and $\varphi_2\wedge\varphi_3$ (both of type $(2, \, 0)$) are $d$-closed but not $d$-exact on $X$, hence they and their conjugates (the latter being of type $(0, \, 2)$) induce non-zero elements in $H^2_{DR}(X, \, \C)$, while $\varphi_1\wedge\varphi_2$ is $d$-exact on $X$ hence it induces the zero class. On the other hand, the $(1, \, 1)$-forms $\varphi_1\wedge\overline{\varphi_1}$, $\varphi_2\wedge\overline{\varphi_2}$, $\varphi_1\wedge\overline{\varphi_2}$ and $\varphi_2\wedge\overline{\varphi_1}$ are all $d$-closed but not $d$-exact. We easily get that $H^2_{DR}(X, \, \C)$ is the $\C$-vector space generated as follows\!\!:

\begin{eqnarray}\label{eqn:H2-Iwa}\nonumber H^2_{DR}(X, \, \C) = &  &\langle\{\varphi_1\wedge\varphi_3\}, \, \{\varphi_2\wedge\varphi_3\}\rangle\\
\nonumber   & \oplus & \langle\{\varphi_1\wedge\overline{\varphi_1}\}, \, \{\varphi_2\wedge\overline{\varphi_2}\}, \, \{\varphi_1\wedge\overline{\varphi_2}\}, \, \{\varphi_2\wedge\overline{\varphi_1}\}\rangle\\
  & \oplus & \langle\{\overline{\varphi_1}\wedge\overline{\varphi_3}\}, \, \{\overline{\varphi_2}\wedge\overline{\varphi_3}\}\}\rangle,  \hspace{3ex} \mbox{hence}\hspace{1ex} b_2(X)=8.\end{eqnarray}

\noindent We have recalled (\ref{eqn:H1-Iwa}) and (\ref{eqn:H2-Iwa}) since they will be used further on.

\vspace{3ex}

 We now review the argument showing that the Iwasawa manifold is balanced. Actually every compact {\it complex parallelisable} manifold will be seen to be balanced. The point of view presented here is that of [AB91a].

\begin{Obs}\label{Obs:n-1_hol-closed} (cf. e.g. [Nak75, Lemma 1.2.] or [AB91a, Remark 3.1.]) Let $X$ be any compact complex (possibly non-K\"ahler) manifold, $\mbox{dim}_{\C}X=n$. Then for every form $u\in C^{\infty}(X, \, \Lambda^{n-1, \, 0}T^{\star}X)$ such that $\bar\partial u=0$, we have $du=0.$ 

 In other words, every holomorphic $(n-1)$-form is $d$-closed on any compact complex manifold of dimension $n$.

\end{Obs}

\noindent {\it Proof.} Let $u\in C^{\infty}(X, \, \Lambda^{n-1, \, 0}T^{\star}X)$ such that $\bar\partial u=0$. Then $du=\partial u$ is of type $(n, \, 0)$ and $d\bar{u}=\overline{du}=\overline{\partial u}$ is of type $(0, \, n)$. We get

\begin{equation}\label{eqn:nn-pos}i^{n^2}du\wedge d\bar{u}\geq 0 \hspace{2ex} \mbox{as an} \,\, (n, \, n)-\mbox{form on} \hspace{2ex} X\end{equation}

\noindent and

\begin{equation}\label{eqn:Stokes}\int\limits_Xi^{n^2}du\wedge d\bar{u} = i^{n^2}\int\limits_Xd(u\wedge d\bar{u})= 0 \hspace{2ex} \mbox{by Stokes}.\end{equation}

\noindent Thus (\ref{eqn:nn-pos}) and (\ref{eqn:Stokes}) yield $i^{n^2}du\wedge d\bar{u}=0$ everywhere on $X$, hence $du=0$ everywhere on $X$. This proves the contention. To justify (\ref{eqn:nn-pos}), write in local holomorphic coordinates $z_1, \dots , z_n$\!\!: \\

$du=f\, dz_1\wedge\dots dz_n, \hspace{2ex} \mbox{hence} \hspace{2ex} i^{n^2}du\wedge d\bar{u}=|f|^2\, idz_1\wedge d\bar{z}_1\dots idz_n\wedge d\bar{z}_n\geq 0.$ \\

\noindent It is clear that $i^{n^2}du\wedge d\bar{u}=0$ iff $f=0$ iff $du=0$.  \hfill $\Box$

\begin{Cor}\label{Cor:1} Let $X$ be a compact complex manifold, $\mbox{dim}_{\C}X=n$. Suppose we have a form $u\in C^{\infty}(X, \, \Lambda^{n-1, \, 0}T^{\star}X)$ such that $\bar\partial u=0$. 

Then the $(n-1, \, n-1)$-form $i^{(n-1)^2}u\wedge\bar{u}$ satisfies

$$i^{(n-1)^2}u\wedge\bar{u}\geq 0  \hspace{2ex} \mbox{and} \hspace{2ex} d\left(i^{(n-1)^2}u\wedge\bar{u}\right)=0 \hspace{2ex} \mbox{on} \hspace{2ex} X.$$

\end{Cor}

\noindent {\it Proof.} The first inequality is checked to hold for any $(n-1, \, 0)$-form $u$ by a trivial calculation. If $\bar\partial u=0$, then $du=0$ by Observation \ref{Obs:n-1_hol-closed}. Then we also have $d\bar{u}=0$ and the second part follows.  \hfill $\Box$

\vspace{2ex}

 Now recall the following standard notion introduced by Wang [Wan54]. A compact complex manifold $X$ is said to be {\bf complex parallelisable} if its holomorphic tangent bundle $T^{1, \, 0}X$ is trivial. This condition is, of course, equivalent to the sheaf of germs of holomorphic $1$-forms $\Omega^1_X$ being trivial. If $n=\mbox{dim}_{\C}X$, the {\it complex parallelisable} condition is equivalent to the existence of $n$ holomorphic vector fields $\theta_1, \dots , \theta_n\in H^0(X, \,T^{1, \, 0}X)$ that are linearly independent at every point of $X$. It is again equivalent to the existence of $n$ holomorphic $1$-forms $\varphi_1, \dots , \varphi_n\in H^0(X, \, \Omega^1_X)$ that are linearly independent at every point of $X$.  

 By a result of Wang [Wan54], every compact {\it complex parallelisable} manifold is the compact quotient $X=G/\Gamma$ of a simply connected, connected {\it complex} Lie group $G$ by a discrete subgroup $\Gamma\subset G$. Conversely, it is obvious that any such quotient is {\it complex parallelisable}. In particular, for any compact {\it complex parallelisable} manifold $X$, $H^0(X, \, T^{1, \, 0}X)\simeq \mathfrak{g}$ where $\mathfrak{g}$ is the Lie algebra of $G$. A compact {\it complex parallelisable} manifold $X$ is said to be {\it nilpotent} (resp. {\it solvable})\footnote{or to be a compact {\it complex parallelisable nilmanifold} (resp. {\it solvmanifold})} if the corresponding complex Lie group $G$ is {\it nilpotent} (resp. {\it solvable}). The Heisenberg group defining the Iwasawa manifold being nilpotent, the Iwasawa manifold is a nilpotent compact {\it complex parallelisable} manifold.

\begin{Cor}\label{Cor:2} (cf. e.g. [AB91a, Remark 3.1.]) Every compact {\bf complex parallelisable} manifold is {\bf balanced}. 

 In particular, the Iwasawa manifold is balanced. 

\end{Cor}

\noindent {\it Proof.} Let $X$ be an arbitrary compact {\it complex parallelisable} manifold, $\mbox{dim}_{\C}X=n$. Let $\varphi_1, \dots , \varphi_n\in H^0(X, \, \Omega^1_X)$ be $n$ holomorphic $1$-forms that are linearly independent at every point of $X$. Consider the $(n-1, \, n-1)$-form on $X$\!\!:

$$\Omega:=i^{(n-1)^2}\sum\limits_{i=1}^n\varphi_1\wedge \dots \wedge \widehat{\varphi_i}\wedge\dots\wedge\varphi_n\wedge\bar{\varphi}_1\wedge \dots \wedge \widehat{\bar{\varphi}_i}\wedge\dots\wedge\bar{\varphi}_n = \sum\limits_{i=1}^n i^{(n-1)^2}u_i\wedge\bar{u}_i,$$

\noindent where $u_i:=\varphi_1\wedge \dots \wedge \widehat{\varphi_i}\wedge\dots\wedge\varphi_n\in C^{\infty}(X, \, \Lambda^{n-1, \, 0}T^{\star}X)$ and \,\, $\widehat{}$\,\, indicates a missing factor. Since $\bar\partial\varphi_k=0$ for all $k=1, \dots , n$, we see that $\bar\partial u_i=0$ for all $i=1, \dots , n$. Then Observation \ref{Obs:n-1_hol-closed} gives $du_i=0$ for all $i=1, \dots , n$, while Corollary \ref{Cor:1} gives

$$\Omega\geq 0   \hspace{2ex}  \mbox{and}  \hspace{2ex} d\Omega=0 \hspace{2ex}  \mbox{on}  \,\, X.$$

\noindent Furthermore, since $\varphi_1, \dots , \varphi_n$ are linearly independent at every point of $X$, we must even have $\Omega > 0.$ Thus $\Omega$ is a $C^{\infty}$ $(n-1, \, n-1)$-form on $X$ satisfying

$$\Omega>0 \hspace{2ex}  \mbox{and} \hspace{2ex} d\Omega=0 \hspace{2ex}  \mbox{on} \,\, X.$$





\noindent Applying Michelsohn's procedure [Mic82] for extracting the $(n-1)^{st}$ root of a smooth positive-definite $(n-1, \, n-1)$-form, there exists a unique $C^{\infty}$ positive-definite $(1, \, 1)$-form $\omega>0$ on $X$ such that $\omega^{n-1}=\Omega$. Since $d(\omega^{n-1})=d\Omega=0$, we see that $\omega$ is a {\it balanced} metric on $X$. The proof is complete. \hfill $\Box$

\vspace{2ex}

Note, however, that very few compact {\it complex parallelisable} manifolds are K\"ahler thanks to a result of Wang.

\begin{Rem}\label{Obs:Wang-K-cor}([Wan54, Corollary 2, p. 776]) Let $X=G/\Gamma$ be a compact complex parallelisable manifold. Then \\

\hspace{3ex} $X$ is K\"ahler \hspace{2ex} iff \hspace{2ex} $G$ is abelian \hspace{2ex} iff \hspace{2ex} $X$ is a complex torus.

\end{Rem}

\vspace{2ex}

\noindent {\bf The Kuranishi family of the Iwasawa manifold} (after Nakamura [Nak75]) \\

\noindent $(a)$\,\, Let $X$ be a compact complex manifold, $\mbox{dim}_{\C}X=n$. Since there are no non-zero $\bar\partial$-exact $(1, \, 0)$-forms on $X$ (for obvious bidegree reasons), we have 

$$H^{1, \, 0}(X, \, \C)=\{u\in C^{\infty}(X, \, \Lambda^{1, \, 0}T^{\star}X) \,\, ; \,\, \bar\partial u=0\},$$

\noindent i.e. $H^{1, \, 0}(X, \, \C)$ consists of holomorphic $1$-forms on $X$. Denoting $h^{1, \, 0}(X):=\mbox{dim}_{\C}H^{1, \, 0}(X, \, \C)$ we have the trivial

\begin{Obs}\label{Obs:trivial-h10} If $X$ is complex parallelisable, then $h^{1, \, 0}(X)=n$.

\end{Obs}  

 \noindent {\it Proof.} By the {\it complex parallelisable} hypothesis on $X$, the rank-$n$ analytic sheaf $\Omega_X^1$ is trivial, hence it is generated by $n$ holomorphic $1$-forms $\varphi_1, \dots , \varphi_n\in H^{1, \, 0}(X, \, \C)$ that are linearly independent at every point of $X$. In particular, $\{\varphi_1, \dots , \varphi_n\}$ is a basis of $H^{1, \, 0}(X, \, \C)\simeq H^0(X, \, \Omega^1_X)$.  \hfill $\Box$

\vspace{2ex}

 Suppose now that $X$ is {\it complex parallelisable}. Let $\theta_1, \dots , \theta_n\in H^0(X, \, T^{1, \, 0}X)$ be holomorphic vector fields that are linearly independent at every point of $X$, chosen to be dual to the holomorphic $(1, \, 0)$-forms $\varphi_1, \dots , \varphi_n\in H^{1, \, 0}(X, \, \C)$ of the above proof. For every smooth function $g:X\to\C$, we have

\begin{equation}\label{eqn:formulae-dd-bar}\partial g=\sum\limits_{\lambda=1}^n(\theta_{\lambda}g)\,\varphi_{\lambda},  \hspace{3ex} \bar\partial g=\sum\limits_{\lambda=1}^n(\bar\theta_{\lambda}g)\,\overline{\varphi}_{\lambda},\end{equation}

\noindent i.e. the familiar formalism induced by local holomorphic coordinates has a global analogue on a compact {\it complex parallelisable} manifold in a formalism where $\theta_{\lambda}$ replaces $\partial/\partial z_{\lambda}$ and $\varphi_{\lambda}$ replaces $dz_{\lambda}$. Thus any $(0, \, 1)$-form $\varphi$ on $X$ has a unique decomposition 

$$\varphi=\sum\limits_{\lambda=1}^nf_{\lambda}\overline{\varphi}_{\lambda}$$

\noindent with $f_1, \dots , f_n: X\to\C$ functions on $X$. Thus there is an implicit $L^2$ inner product on $C^{\infty}(X, \, \Lambda^{0, \, 1}T^{\star}X)$ defined as follows (no Hermitian metric is needed on $X$): for any $\varphi=\sum\limits_{\lambda=1}^nf_{\lambda}\overline{\varphi}_{\lambda}, \psi=\sum\limits_{\lambda=1}^ng_{\lambda}\overline{\varphi}_{\lambda}\in C^{\infty}(X, \, \Lambda^{0, \, 1}T^{\star}X),$ set 

\begin{equation}\label{eqn:L2innerproduct}\langle\langle\varphi, \, \psi\rangle\rangle :=\int\limits_X\bigg(\sum\limits_{\lambda=1}^nf_{\lambda}\,\bar{g}_{\lambda}\bigg)\, i^{n^2}\varphi_1\wedge\dots\wedge\varphi_n\wedge\overline{\varphi}_1\wedge\dots\wedge\overline{\varphi}_n.\end{equation}

\noindent It is clear that $dV:=i^{n^2}\varphi_1\wedge\dots\wedge\varphi_n\wedge\overline{\varphi}_1\wedge\dots\wedge\overline{\varphi}_n>0$ is a global volume form on $X$ and that the above $L^2$ inner product is independent of the choices made. We can define the formal adjoint $\bar\partial^{\star}$ of $\bar\partial$ w.r.t. this $L^2$ inner product in the usual way: for any smooth $(0, \, 1)$-form $\varphi$, define $\bar\partial^{\star}\varphi$ to be the unique smooth function on $X$ satisfying

$$\langle\langle\bar\partial^{\star}\varphi, \, g\rangle\rangle = \langle\langle\varphi, \, \bar\partial g\rangle\rangle$$

\noindent for any smooth function $g$ on $X$. A trivial calculation using Stokes's theorem gives

\begin{equation}\label{eqn:d-baradjoint}\bar\partial^{\star}\varphi = - \sum\limits_{\lambda=1}^n\theta_{\lambda}f_{\lambda}\end{equation}

\noindent for any smooth $(0, \, 1)$-form $\varphi=\sum\limits_{\lambda=1}^nf_{\lambda}\overline{\varphi}_{\lambda}$ on $X$. Thus we see that

\begin{equation}\label{eqn:d-barstar-calc}\bar\partial^{\star}\overline{\varphi}_{\nu}=0, \hspace{3ex} \nu=1, \dots , n,\end{equation}

\noindent because $\overline{\varphi}_{\nu}=\sum\limits_{\lambda=1}^n\delta_{\nu\lambda}\overline{\varphi}_{\lambda}$ and $\theta_{\lambda}\delta_{\nu\lambda}=0$ (since the $\delta_{\nu\lambda}$ are constants). 

 Now denote by $r\in\{0, 1, \dots , n\}$ the number of $d$-closed forms among $\varphi_1, \dots , \varphi_n$. After a possible reordering, we can suppose that $\varphi_1, \dots , \varphi_r$ are $d$-closed and $\varphi_{r+1}, \dots, \varphi_n$ are not $d$-closed. Then we have

\begin{equation}\label{eqn:d-bar-closed}\partial\varphi_1= \dots = \partial\varphi_r =0 \hspace{2ex} \mbox{or equivalently} \hspace{2ex} \bar\partial\overline{\varphi}_1= \dots = \bar\partial\overline{\varphi}_r =0.\end{equation}

\noindent Thus the $\bar\partial$-closed $(0, \, 1)$-forms $\overline{\varphi}_1, \dots , \overline{\varphi}_r$ define Dolbeault $(0, \, 1)$-cohomology classes in $H^{0, \, 1}(X, \, \C)$. 

 We can define the $\bar\partial$-Laplacian on forms of $X$ in the usual way\!\!:

$$\Delta'':=\bar\partial\bar\partial^{\star} + \bar\partial^{\star}\bar\partial.$$

\noindent The corresponding harmonic space of $(0, \, 1)$-forms ${\cal H}^{0, \, 1}_{\Delta''}(X, \, \C):=\ker\Delta''=\ker\bar\partial\cap\ker\bar\partial^{\star}$ satisfies the Hodge isomorphism ${\cal H}^{0, \, 1}_{\Delta''}(X, \, \C)\simeq H^{0, \, 1}(X, \, \C)$. Notice that (\ref{eqn:d-barstar-calc}) and (\ref{eqn:d-bar-closed}) give

\begin{equation}\label{eqn:harm}\Delta''\overline{\varphi}_{\nu}=0, \hspace{3ex} \nu=1, \dots , r,\end{equation}

\noindent i.e. the forms $\overline{\varphi}_1, \dots , \overline{\varphi}_r$ are $\Delta''$-harmonic. On the other hand, $\overline{\varphi}_{r+1}, \dots , \overline{\varphi}_n$ are not $\Delta''$-harmonic. Thus the number $r$ of linearly independent $d$-closed holomorphic $1$-forms of $X$ (independent of the choice of $\varphi_1, \dots , \varphi_n$) satisfies:

\begin{equation}\label{eqn:r-h01-ineq}r\leq h^{0, \, 1}(X).\end{equation}

\vspace{1ex}

 Suppose now that the compact {\it complex parallelisable} $X$ is {\it nilpotent}.

\begin{Fact}\label{Fact:Chev-decomp}(see e.g. [Nak75] or [CFGU00, p.5405-5406])
If $X$ is a compact {\bf complex parallelisable nilmanifold}, the holomorphic $1$-forms $\varphi_1, \dots , \varphi_n$ that form a basis of $H^{1, \, 0}(X, \, \C)$ can be chosen such that

\begin{equation}\label{eqn:Chev-decomp}d\varphi_{\mu}=\sum\limits_{1\leq\lambda<\nu\leq n}c_{\mu\lambda\nu}\,\varphi_{\lambda}\wedge\varphi_{\nu},\hspace{3ex} 1\leq\mu\leq n,\end{equation}

\noindent with constant coefficients $c_{\mu\lambda\nu}\in\C$ satisfying

$$c_{\mu\lambda\nu}=0  \hspace{2ex} \mbox{whenever} \hspace{2ex} \mu\leq\lambda \hspace{1ex} \mbox{or} \hspace{1ex} \mu\leq\nu.$$

\end{Fact}

 Taking this standard fact (which in [Nak75] follows from the existence of a {\it Chevalley decomposition} of the nilpotent Lie algebra $\mathfrak{g}$) for granted, we now spell out the details of the proof of the following result of Kodaira along the lines of [Nak75, Theorem 3, p. 100].

\begin{The}\label{The:Kodaira}(Kodaira) If $X$ is a compact {\bf complex parallelisable nilmanifold}, then $h^{0, \, 1}(X)=r$.

 Moreover, the $\Delta''$-harmonic $(0, \, 1)$-forms $\overline{\varphi}_1,\dots , \overline{\varphi}_r$ form a basis of the harmonic space ${\cal H}_{\Delta''}^{0, \, 1}(X, \, \C)$. Equivalently, the Dolbeault $(0, \, 1)$-cohomology classes $\{\overline{\varphi}_1\}, \dots , \{\overline{\varphi}_r\}$ form a basis of $H^{0, \, 1}(X, \, \C)$. 

\end{The}

\noindent {\it Proof.} The only thing that needs proving is that the linearly independent forms $\overline{\varphi}_1,\dots , \overline{\varphi}_r\in{\cal H}_{\Delta''}^{0, \, 1}(X, \, \C)$ generate ${\cal H}_{\Delta''}^{0, \, 1}(X, \, \C)$. Pick an arbitrary $C^{\infty}$ $(0, \, 1)$-form $\varphi$ on $X$ and write

$$\varphi=\sum\limits_{\lambda=1}^nf_{\lambda}\, \overline{\varphi}_{\lambda}$$

\noindent with $C^{\infty}$ functions $f_1, \dots , f_n$ on $X$. Using formula (\ref{eqn:formulae-dd-bar}) for $\bar\partial$ and the obvious identities $\bar\partial\,\overline{\varphi}_{\lambda} = d\,\overline{\varphi}_{\lambda}$, $\lambda=1, \dots , n$, due to ${\varphi}_{\lambda}$ being holomorphic, we get\!\!:

\begin{eqnarray}\nonumber\bar\partial\varphi & = & \sum\limits_{\lambda, \,\nu=1}^n(\bar\theta_{\nu}f_{\lambda})\, \overline{\varphi}_{\nu}\wedge\overline{\varphi}_{\lambda} + \sum\limits_{\mu=1}^nf_{\mu}\, d\,\overline{\varphi}_{\mu} \\
\nonumber & = & \sum\limits_{\lambda, \,\nu=1}^n(\bar\theta_{\nu}f_{\lambda})\, \overline{\varphi}_{\nu}\wedge\overline{\varphi}_{\lambda} + \sum\limits_{\mu=1}^nf_{\mu}\sum\limits_{1\leq\nu<\lambda\leq n}\overline{c_{\mu\nu\lambda}}\, \overline{\varphi}_{\nu}\wedge\overline{\varphi}_{\lambda}\\
\label{eqn:d-bar-phi}  & = & \sum\limits_{1\leq\nu < \lambda\leq n}\bigg(\bar{\theta}_{\nu}f_{\lambda}-\bar{\theta}_{\lambda}f_{\nu} + \sum\limits_{\mu=1}^n\overline{c_{\mu\nu\lambda}}\,f_{\mu}\bigg)\overline{\varphi}_{\nu}\wedge\overline{\varphi}_{\lambda},\end{eqnarray}

\noindent where the second line above follows from the conjugate of (\ref{eqn:Chev-decomp}). 

 Now $\varphi$ is $\Delta''$-harmonic if and only if \\

\noindent $(i)$\,\, $\bar\partial\varphi=0 \Longleftrightarrow \bar{\theta}_{\nu}f_{\lambda}-\bar{\theta}_{\lambda}f_{\nu} + \sum\limits_{\mu=1}^n\overline{c_{\mu\nu\lambda}}\,f_{\mu}=0$ \hspace{2ex} for $1\leq\nu<\lambda\leq n$  (cf. (\ref{eqn:d-bar-phi})); \\

\noindent and \\

\noindent $(ii)$\,\, $\bar\partial^{\star}\varphi=0 \Longleftrightarrow \sum\limits_{\lambda=1}^n\theta_{\lambda}\, f_{\lambda}=0$  \hspace{2ex}  (cf. (\ref{eqn:d-baradjoint})).

 \vspace{1ex}

 Suppose that $\varphi$ is $\Delta''$-harmonic. Then the above $(i)$ reads\!\!:

$$\bar{\theta}_{\lambda}f_{\nu} = \sum\limits_{\mu=1}^n\overline{c_{\mu\nu\lambda}}\,f_{\mu} + \bar{\theta}_{\nu}f_{\lambda}, \hspace{2ex} 1\leq\nu<\lambda\leq n.$$

\noindent Summing over $\lambda=1, \dots , n$ and using formula (\ref{eqn:formulae-dd-bar}) for $\bar\partial$, we get

$$\bar\partial f_{\nu}=\sum\limits_{\lambda=1}^n(\bar{\theta}_{\lambda}f_{\nu})\, \overline{\varphi}_{\lambda} = \sum\limits_{\lambda, \,\mu=1}^n\overline{c_{\mu\nu\lambda}}\,f_{\mu}\overline{\varphi}_{\lambda} + \sum\limits_{\lambda=1}^n(\bar{\theta}_{\nu}f_{\lambda})\,\overline{\varphi}_{\lambda}, \hspace{2ex} \nu=1, \dots , n,$$

\noindent with the understanding that $c_{\mu\nu\lambda}=0$ if $\nu\geq\lambda$. Now $\Delta''f_{\nu}=\bar\partial^{\star}\bar\partial f_{\nu}$ since $f_{\nu}$ is a function. Taking $\bar\partial^{\star}$ on either side above and using formula (\ref{eqn:d-baradjoint}) for $\bar\partial^{\star}$, we get

\begin{eqnarray}\nonumber\Delta''f_{\nu} & = & -\sum\limits_{\lambda,\,\mu=1}^n\theta_{\lambda}\,(\overline{c_{\mu\nu\lambda}}\, f_{\mu}) - \sum\limits_{\lambda=1}^n\theta_{\lambda}\,(\bar{\theta}_{\nu}f_{\lambda})  \\
\label{eqn:delta''f-nu}  & = & -\sum\limits_{\lambda,\,\mu=1}^n\overline{c_{\mu\nu\lambda}}\, \theta_{\lambda}f_{\mu}, \hspace{2ex} \mbox{for all} \hspace{1ex} \nu=1, \dots , n,\end{eqnarray}

\noindent because $\theta_{\lambda}\,(\overline{c_{\mu\nu\lambda}}f_{\mu}) = \overline{c_{\mu\nu\lambda}}\,\theta_{\lambda}f_{\mu}$ due to $\overline{c_{\mu\nu\lambda}}$ being constant, while $\sum\limits_{\lambda=1}^n\theta_{\lambda}f_{\lambda}=0$ due to $\varphi$ being $\Delta''$-harmonic (cf. $(ii)$ or (\ref{eqn:d-baradjoint})).

 Taking now $\nu=n$ in (\ref{eqn:delta''f-nu}), we get $\Delta''f_n=0$ since $
c_{\mu n\lambda}=0$ for all $\mu, \lambda$ by Fact \ref{Fact:Chev-decomp} and the obvious inequality $\mu\leq\nu=n$. Thus the compactness of $X$ and the ellipticity of $\Delta''$ yield

\begin{equation}\label{eqn:const-n}f_n \hspace{2ex} \mbox{is constant on} \hspace{1ex} X \hspace{1ex} \mbox{if} \hspace{1ex} \Delta''\varphi=0.\end{equation}

 Taking now $\nu=n-1$ in (\ref{eqn:delta''f-nu}) and using the fact that $\theta_{\lambda}f_n=0$ for all $\lambda$ (since $f_n$ is constant by (\ref{eqn:const-n})), we get

$$\Delta''f_{n-1} = -\sum\limits_{\lambda=1}^n\bigg(\sum\limits_{\mu=1}^{n-1}\overline{c_{\mu n-1\lambda}}\, \theta_{\lambda}f_{\mu}\bigg)=0 \hspace{2ex} \mbox{on}\hspace{1ex} X,$$

\noindent since $c_{\mu n-1\lambda}=0$ for all $\mu=1, \dots , n-1$ and $\lambda=1, \dots , n$ by Fact \ref{Fact:Chev-decomp} and the obvious inequality $\mu\leq \nu=n-1$. Hence we get

\begin{equation}\label{eqn:const-n-1}f_{n-1} \hspace{2ex} \mbox{is constant on} \hspace{1ex} X \hspace{1ex} \mbox{if} \hspace{1ex} \Delta''\varphi=0.\end{equation}

 We can now continue by decreasing induction on $\nu$. Taking $\nu=n-2$ in (\ref{eqn:delta''f-nu}) and using the fact that $\theta_{\lambda}f_n = \theta_{\lambda}f_{n-1} = 0$ for all $\lambda$ (since $f_n$ is constant by (\ref{eqn:const-n}) and $f_{n-1}$ is constant by (\ref{eqn:const-n-1})), we get

$$\Delta''f_{n-2} = -\sum\limits_{\lambda=1}^n\bigg(\sum\limits_{\mu=1}^{n-2}\overline{c_{\mu n-2\lambda}}\, \theta_{\lambda}f_{\mu}\bigg)=0 \hspace{2ex} \mbox{on}\hspace{1ex} X,$$

\noindent since $c_{\mu n-2\lambda}=0$ for all $\mu=1, \dots , n-2$ and $\lambda=1, \dots , n$ by Fact \ref{Fact:Chev-decomp} and the obvious inequality $\mu\leq \nu=n-2$. Hence we get

\begin{equation}\label{eqn:const-n-2}f_{n-2} \hspace{2ex} \mbox{is constant on} \hspace{1ex} X \hspace{1ex} \mbox{if} \hspace{1ex} \Delta''\varphi=0.\end{equation}

 Running a decreasing induction on $\nu$, we get

\begin{equation}\label{eqn:const-all}f_{\nu}:= C_{\nu} \hspace{2ex} \mbox{is constant on} \hspace{1ex} X \hspace{1ex} \mbox{for all} \hspace{1ex} \nu=1, \dots , n \hspace{1ex} \mbox{if} \hspace{1ex} \Delta''\varphi=0.\end{equation}

 We conclude that whenever $\Delta''\varphi=0$ we have

$$\varphi=\sum\limits_{\nu=1}^n C_{\nu}\,\overline{\varphi}_{\nu} \hspace{2ex} \mbox{with} \hspace{1ex} C_{\nu} \hspace{1ex} \mbox{constant for all} \hspace{1ex} \nu=1, \dots , n.$$

 On the other hand, since $\Delta''\varphi=0$, we must have $\bar\partial\varphi=0$ which amounts to $\sum\limits_{\nu=1}^n C_{\nu}\,\bar\partial\overline{\varphi}_{\nu}=0$. However, we know that $\bar\partial\overline{\varphi}_{\nu}=0$ for all $\nu\in\{1, \dots , r\}$ (cf. (\ref{eqn:d-bar-closed})), hence $\sum\limits_{\nu=r+1}^n C_{\nu}\,\bar\partial\overline{\varphi}_{\nu}=0$. Now the forms

$$\bar\partial\overline{\varphi}_{\nu} = d\overline{\varphi}_{\nu}=\sum\limits_{\lambda, \mu}\overline{c_{\nu\lambda\mu}}\,\overline{\varphi}_{\lambda}\wedge\overline{\varphi}_{\mu}, \hspace{2ex} \nu=1, \dots , n,$$

\noindent are linearly independent because $\overline{\varphi}_1, \dots , \overline{\varphi}_n$ are linearly independent at every point of $X$. Hence $C_{\nu}=0$ for all $\nu=r+1, \dots , n$. We get

$$\varphi=\sum\limits_{\nu=1}^r C_{\nu}\,\overline{\varphi}_{\nu} \hspace{2ex} \mbox{with} \hspace{1ex} C_{\nu} \hspace{1ex} \mbox{constant for all} \hspace{1ex} \nu=1, \dots , r.$$

 Since $\varphi$ has been chosen arbitrary in ${\cal H}_{\Delta''}^{0, \, 1}(X, \, \C)$, we have proved that the linearly independent forms $\overline{\varphi}_1,\dots , \overline{\varphi}_r\in{\cal H}_{\Delta''}^{0, \, 1}(X, \, \C)$ generate ${\cal H}_{\Delta''}^{0, \, 1}(X, \, \C)$. The proof of Kodaira's theorem \ref{The:Kodaira} is complete.  \hfill $\Box$

\vspace{3ex}

 When applying Observation \ref{Obs:trivial-h10} and Kodaira's Theorem \ref{The:Kodaira} to the Iwasawa manifold (a compact {\it complex parallelisable nilmanifold} of dimension $n=3$ having $r=2$), we get the following classical fact.

\begin{Obs}\label{Obs:coh-Iwasawa} For the Iwasawa manifold, we have

$$h^{1, \, 0}= 3  \hspace{3ex} \mbox{and} \hspace{3ex} h^{0, \, 1}=2.$$

 Since, on the other hand, the first Betti number is $b_1=4$ (see e.g. (\ref{eqn:H1-Iwa})), we see that $b_1<h^{1, \, 0} + h^{0, \, 1}$. Thus the Fr\"olicher spectral sequence of the Iwasawa manifold does not degenerate at $E_1$\footnote{Historically, the Iwasawa manifold was the first example exhibited of a compact complex manifold with this property. It is even known to have $E_1\neq E_2=E_{\infty}$ in its Fr\"olicher spectral sequence.}. In particular, the $\partial\bar\partial$-lemma does not hold on the Iwasawa manifold.

\end{Obs}

 On the other hand, we have seen in Corollary \ref{Cor:2} that the Iwasawa manifold is balanced, hence sG. Combined with Observation \ref{Obs:coh-Iwasawa}, this fact proves part $(a)$ of Theorem \ref{The:sG-specseq-unrel}. The same statement will be proved again in section \ref{section:ES} (see comments after Gauduchon's Theorem \ref {The:Gau-twist}) by a twistor space featuring as the central fibre in the Eastwood-Singer family [ES93].

\vspace{3ex}

\noindent $(b)$\,\,  Now suppose that $X$ is the {\it Iwasawa manifold}. In particular, $X$ is a compact {\it complex parallelisable nilmanifold} of complex dimension $3$. Let $\varphi_1=dz_1$, $\varphi_2=dz_2$, $\varphi_3=dz_3 - z_1dz_2$ be the holomrophic $1$-forms on $X$ defined in (\ref{eqn:Iwa-forms}); they are linearly independent at every point of $X$. Since $\varphi_1$ and $\varphi_2$ are $d$-closed while $\varphi_3$ is not $d$-closed, $r=2$ for the Iwasawa manifold. By Kodaira's theorem \ref{The:Kodaira}, the $\C$-vector space $H^{0, \, 1}(X, \, \C)$ has complex dimension $2$ and is spanned by the Dolbeault cohomology classes $\{\overline{\varphi}_1\}$ and $\{\overline{\varphi}_2\}$. Let $\theta_1, \theta_2, \theta_3\in H^0(X, \, \Omega^1_X)$ be the holomorphic vector fields dual to $\varphi_1, \varphi_2, \varphi_3$. They are given by

\begin{equation}\label{eqn:Iwa-hol-vectfiel}\theta_1=\frac{\partial}{\partial z_1}, \hspace{1ex} \theta_2=\frac{\partial}{\partial z_2} + z_1\frac{\partial}{\partial z_3}, \hspace{1ex} \theta_3=\frac{\partial}{\partial z_3}\end{equation}

\noindent and satisfy the relations

\begin{equation}\label{eqn:hol-vectfiel-comm}[\theta_1, \, \theta_2] = - [\theta_2, \, \theta_1] = \theta_3, \hspace{3ex} [\theta_1, \, \theta_3] = [\theta_2, \, \theta_3] = 0,\end{equation}

\noindent i.e. $[\theta_i, \, \theta_j]=0$ whenever $\{i, \, j\}\neq \{1, \, 2\}$.

 Since the holomorphic tangent bundle $T^{1, \, 0}X$ is trivial and spanned by $\theta_1, \theta_2, \theta_3$, the cohomology group $H^{0, \, 1}(X, \, T^{1, \, 0}X)$ of $T^{1, \, 0}X$-valued $(0, \, 1)$-forms on $X$ is a $\C$-vector space of dimension $6$ spanned by the classes of $\theta_i\,\overline{\varphi}_{\lambda}:$

\begin{equation}\label{eqn:H01T10}H^{0, \, 1}(X, \, T^{1, \, 0}X)=\bigoplus\limits_{1\leq i \leq 3, \, 1\leq \lambda\leq 2}\C\{\theta_i\,\overline{\varphi}_{\lambda}\}, \hspace{3ex} \mbox{dim}_{\C}H^{0, \, 1}(X, \, T^{1, \, 0}X)=6.\end{equation}

 This will be seen to imply that the Kuranishi family of the Iwasawa manifold is a $6$\,-parameter family.

\vspace{3ex}

$(c)$\,\,  We now briefly recall a few basic facts in the deformation theory of compact complex manifolds in order to fix the notation. The standard reference is Kodaira's book [Kod86]. Given a compact complex manifold $X$ of complex dimension $n$, let $\bar\partial$ be the Cauchy-Riemann operator representing the complex structure $J=J_0$ of $X$ and let $z_1, \dots , z_n$ be local holomorphic coordinates about an arbitrary point on $X$. A deformation $X_t$ of the complex structure of $X$ is represented by a vector $(0, \, 1)$-form 

$$\psi(t)\in C^{\infty}(X, \, \Lambda^{0, \, 1}T^{\star}X\otimes T^{1, \, 0}X)$$

\noindent in the following sense: the Cauchy-Riemann operator $\bar\partial_t$ representing the (almost) complex structure $J_t$ of $X_t$ is defined by

$$\bar\partial_t:=\bar\partial - \psi(t).$$

\noindent Equivalently, a locally defined $C^{\infty}$ function $f$ on $X$ satisfies the equivalence\!\!:

$$f \,\, \mbox{is} \,\, J_t\!-\!\mbox{holomorphic}  \hspace{2ex} \mbox{if and only if} \hspace{2ex} (\bar\partial - \psi(t))f=0.$$

\noindent The integrability condition reads\!\!:

\begin{eqnarray}\label{eqn:int-cond}J_t \,\, \mbox{is integrable} \hspace{1ex} & \mbox{if and only if} & \hspace{1ex} \bar\partial\psi(t)=\frac{1}{2}[\psi(t), \, \psi(t)] \\
\nonumber  & \mbox{if and only if} & \hspace{1ex} \mbox{the system of n PDE's} \hspace{1ex} \bar\partial_tf=0 \hspace{2ex} \mbox{has}\\
\nonumber  &  & \hspace{1ex} n\,\, \mbox{linearly independent}\,\, C^{\infty} \,\, \mbox{solutions in}\\
 \label{eqn:int-cond2}     &  & \hspace{1ex} \mbox{a neighbourhood of any point of}\,\, X.\end{eqnarray}

 Recall furthermore that the {\it Kodaira-Spencer map} of a holomorphic family $\pi : {\cal X} \longrightarrow B\subset\C^m$ at $t=0$ (having supposed that $0\in B$) is

$$\rho_0:T_0B\longrightarrow H^{0, \, 1}(X, \, T^{1, \, 0}X),  \hspace{3ex} \frac{\partial}{\partial t}\mapsto -\frac{\partial\psi(t)}{\partial t}_{|t=0}:=\frac{\partial X_t}{\partial t}_{|t=0}.$$

\noindent Thus $\frac{\partial X_t}{\partial t}_{|t=0}$ denotes the infinitesimal deformation of $X$ at $t=0$. Recall now the fundamental Kuranishi Theorem of Existence.

\begin{The}\label{The:Kur}(Kuranishi [Kur62]) Given any compact complex manifold $X$, there exists a {\bf complete} holomorphic family

$$\pi : {\cal X}\longrightarrow B\subset\Delta_{\varepsilon}:=\{t\in\C^m\,\, ; \,\, |t|<\varepsilon\}$$

\noindent such that $X_0=X$, for some small $\varepsilon >0$ and some {\bf analytic} subset $B$ of the ball $\Delta_{\varepsilon}$, where $m:=\mbox{dim}_{\C}H^{0, \, 1}(X, \, T^{1, \, 0}X)$.

\end{The}

 {\it Complete} means that the Kodaira-Spencer map $\rho_0:T_0B\longrightarrow H^{0, \, 1}(X, \, T^{1, \, 0}X)$ is surjective. This family is called the {\it Kuranishi family} of $X$. The base $B$ of the family may have singularities and arises as

$$B=\{t\in\Delta_{\varepsilon}\,\, ; \,\, f_1(t)= \dots = f_l(t)=0\},$$

\noindent where $l:=\mbox{dim}_{\C}H^{0, \, 2}(X, \, T^{1, \, 0}X)$ and $f_1, \dots , f_l:\Delta_{\varepsilon}\to\C$ are holomorphic functions. In the special case of a manifold $X$ satisfying $H^{0, \, 2}(X, \, T^{1, \, 0}X)=0$, the base $B$ is smooth and Kuranishi's theorem reduces to the earlier Kodaira-Nirenberg-Spencer theorem of existence. 

 The construction of the Kuranishi family of a given $X$ amounts to the construction of a family of vector $(0, \, 1)$-forms $\psi(t)\in C^{\infty}(X, \, \Lambda^{0, \, 1}T^{\star}X\otimes T^{1, \, 0}X)$ satisfying the integrability condition (\ref{eqn:int-cond}) for $t=(t_1, \dots , t_m)$ in the largest possible subset of some $\Delta_{\varepsilon}\subset\C^m$. Given an arbitrary basis $\{\beta_1, \dots , \beta_m\}$ of $H^{0, \, 1}(X, \, T^{1, \, 0}X)$, set 

\begin{equation}\label{eqn:psi1-basis}\psi_1(t):=t_1\beta_1 + \dots t_m\beta_m\in H^{0, \, 1}(X, \, T^{1, \, 0}X).\end{equation} 

\noindent Identifying $H^{0, \, 1}(X, \, T^{1, \, 0}X)$ with ${\cal H}^{0, \, 1}_{\Delta''}(X, \, T^{1, \, 0}X)$ by the Hodge isomorphism, we see that $\Delta''\psi_1(t)=0$ (i.e. $\psi_1(t)$ is $\Delta''$-harmonic). 

 Since $\psi_1$ need not satisfy the integrability condition (\ref{eqn:int-cond}), we search for a power series

\begin{equation}\label{eqn:power-series}\psi(t)=\psi_1(t) + \sum\limits_{\nu=2}^{+\infty}\psi_{\nu}(t),\end{equation}

\noindent where $\psi_{\nu}(t)=\sum\limits_{\nu_1 + \dots + \nu_m=\nu}\psi_{\nu_1\dots\nu_m}t_1^{\nu_1}\dots t_m^{\nu_m}$ is a homogeneous polynomial of degree $\nu$ in $t_1, \dots , t_m$ whose coefficients $\psi_{\nu_1\dots\nu_m}\in C^{\infty}(X, \, \Lambda^{0, \, 1}T^{\star}X\otimes T^{1, \, 0}X)$ will be determined such that the following two conditions are fulfilled: \\

\noindent $\bullet$ the power series defining $\psi(t)$ converges and its sum is $C^{\infty}$ on $X\times\Delta_{\varepsilon}$ for some small $\varepsilon >0$;

 Kuranishi's proof achieves convergence in a H\"older norm $|\,\,\,\,|_{k, \,\alpha}$ for all $k\geq 2$ and all $t\in\Delta_{\varepsilon}$ provided that $\varepsilon >0$ is small enough. 

\vspace{1ex}

\noindent $\bullet$ the integrability condition $\bar\partial\psi(t)=\frac{1}{2}[\psi(t), \, \psi(t)]$ holds (cf. (\ref{eqn:int-cond})).

\vspace{1ex}

 To fulfill the integrability condition (\ref{eqn:int-cond}), it suffices to ensure that $[\psi(t), \, \psi(t)]\in C^{\infty}(X, \, \Lambda^{0, \, 2}T^{\star}X\otimes T^{1, \, 0}X)$ is $\bar\partial$-exact and that $\psi(t)-\psi_1(t)$ is the minimal $L^2$-norm solution of equation $\bar\partial u=\frac{1}{2}[\psi(t), \, \psi(t)]$ (recall that $\bar\partial\psi_1(t)=0$). Minimality of the solution's $L^2$-norm translates to

$$\psi(t)=\psi_1(t) + \frac{1}{2}\bar\partial^{\star}\Delta^{''-1}[\psi(t), \, \psi(t)],$$

\noindent a formula that is easily seen to be equivalent to

\begin{equation}\label{eqn:psi-nu-form}\psi_{\nu}(t)= \frac{1}{2}\sum\limits_{\mu=1}^{\nu -1}\bar\partial^{\star}\Delta^{''-1}[\psi_{\mu}(t), \, \psi_{\nu - \mu}(t)] \hspace{2ex} \mbox{for all}\,\,\nu\geq 2.\end{equation}

\noindent This means that $\psi_{\nu}(t)$ is the minimal $L^2$-norm solution of the equation $\bar\partial u=v_{\nu}$, where $v_{\nu}$ is the projection of $1/2\sum_{1\leq\mu\leq\nu-1}[\psi_{\mu}(t), \, \psi_{\nu - \mu}(t)]$ onto $\mbox{Im}\,\bar\partial$. In particular, $\psi_{\nu}(t)\in\mbox{Im}\,\bar\partial^{\star}$ for all $\nu\geq 2$. 

 Identities (\ref{eqn:psi-nu-form}) allow one to construct $\psi_{\nu}(t)$, $\nu\geq 2$, inductively from $\psi_1(t)$ defined in (\ref{eqn:psi1-basis}). Convergence in H\"older norm $|\,\,\,\,|_{k, _, \alpha}$ of the resulting series (\ref{eqn:power-series}) follows from {\it a priori estimates} on the Laplacian $\Delta''$, while the integrability condition (\ref{eqn:int-cond}) for the sum $\psi(t)$ of this series is seen to be equivalent to

\begin{equation}\label{eqn:intcond-harmproj}H[\psi(t), \, \psi(t)]=0,\end{equation}

\noindent where $H: C^{\infty}(X, \, \Lambda^{0, \, 2}T^{\star}X\otimes T^{1, \, 0}X)\to{\cal H}^{0, \, 2}_{\Delta''}(X, \, T^{1, \, 0}X)$ is the harmonic projector. Condition (\ref{eqn:intcond-harmproj}) requires $[\psi(t), \, \psi(t)]$ to have no harmonic component which, for a $\bar\partial$-closed form, is equivalent to $\bar\partial$-exactness (precisely what is needed in view of (\ref{eqn:int-cond})). If $\{\gamma_1, \dots , \gamma_l\}$ is any orthonormal basis of ${\cal H}^{0, \, 2}_{\Delta''}(X, \, T^{1, \, 0}X)$, then

$$H[\psi(t), \, \psi(t)]=\sum\limits_{k=1}^l\langle[\psi(t), \, \psi(t)], \, \gamma_k\rangle\, \gamma_k, \hspace{2ex} t\in\Delta_{\varepsilon},$$

\noindent and we see that the vanishing condition (\ref{eqn:intcond-harmproj}) is equivalent to $f_1(t)=\dots = f_l(t)=0$, where $f_k(t):=\langle[\psi(t), \, \psi(t)], \, \gamma_k\rangle$ for all $k=1, \dots , l$ and $t\in\Delta_{\varepsilon}$. Thus the integrability condition (\ref{eqn:int-cond}) is satisfied for $t\in B$, where

\begin{equation}\label{eqn:BvsDelta}B:=\{t\in\Delta_{\varepsilon}\,\, ; \,\, f_1(t)= \dots = f_l(t)=0\}\subset\Delta_{\varepsilon}\end{equation}

\noindent is analytic.

\vspace{3ex}

\noindent $(d)$\,\, {\it Nakamura's calculation of the Kuranishi family of the Iwasawa manifold} \\

 Given a compact complex manifold $X$ of dimension $n$, for any vector $(0, \, 1)$-forms $\psi, \tau\in C^{\infty}(X, \, \Lambda^{0, \, 1}T^{\star}X\otimes T^{1, \, 0}X)$ written locally as

$$\psi=\sum\limits_{\alpha=1}^n\psi^{\alpha}\frac{\partial}{\partial z_{\alpha}}, \hspace{2ex} \tau=\sum\limits_{\beta=1}^n\tau^{\beta}\frac{\partial}{\partial z_{\beta}}, \hspace{3ex} \psi^{\alpha}, \tau^{\beta}\in C^{\infty}(X, \, T^{1, \, 0}X),$$

\noindent Kuranishi defines in general

$$[\psi, \, \tau]:=\sum\limits_{\alpha, \,\beta=1}^n\bigg(\psi^{\alpha}\wedge\frac{\partial\tau^{\beta}}{\partial z_{\alpha}} + \tau^{\alpha}\wedge\frac{\partial\psi^{\beta}}{\partial z_{\alpha}}\bigg)\frac{\partial}{\partial z_{\beta}}.$$

 Now fix $X=\C^3/\Gamma$ to be the Iwasawa manifold. Then $n=3$ and we get\!\!:

\begin{equation}\label{eqn:vector-brackets}[\theta_i\overline{\varphi}_{\lambda}, \, \theta_k\overline{\varphi}_{\nu}] = [\theta_i, \,  \theta_k]\, \overline{\varphi}_{\lambda}\wedge\overline{\varphi}_{\nu}, \hspace{3ex} i, k, \lambda, \nu=1,2,3,\end{equation}

\noindent with $[\theta_i, \,  \theta_k]$ given in (\ref{eqn:hol-vectfiel-comm}). 

 We have seen in (\ref{eqn:H01T10}) that the classes $\{\theta_i\,\overline{\varphi_{\lambda}}\}$, with $1\leq i\leq 3,\, 1\leq\lambda\leq 2$, form a basis of $H^{0, \, 1}(X, \, T^{1, \, 0}X)$. Consequently the Kuranishi family of $X$ can be described by $6$ parameters $t=(t_{i\lambda})_{1\leq i\leq 3, \,\, 1\leq\lambda\leq 2}$. By (\ref{eqn:harm}), the $T^{1, \, 0}X$-valued $(0, \, 1)$-forms $\theta_i\,\overline{\varphi_{\lambda}}$ are $\Delta''$-harmonic when $1\leq\lambda\leq 2$. In order to construct the vector $(0, \, 1)$-forms $\psi(t)\in C^{\infty}(X, \, \Lambda^{0, \, 1}T^{\star}X\otimes T^{1, \, 0}X)$ that describe the Kuranishi family of $X=\C^3/\Gamma$, formula (\ref{eqn:psi1-basis}) prescribes to start off by setting 

\begin{equation}\label{eqn:psi1-def}\psi_1(t):=\sum\limits_{i=1}^3\sum\limits_{\lambda=1}^2t_{i\lambda}\theta_i\overline{\varphi}_{\lambda}, \hspace{2ex} t=(t_{i\lambda})_{1\leq i\leq 3, \,\, 1\leq\lambda\leq 2},\end{equation}

\noindent for which we see that

\begin{equation}\nonumber\label{eqn:psi_1_brackets1}\frac{1}{2}[\psi_1(t), \, \psi_1(t)] =  \frac{1}{2} \sum\limits_{i, j=1, 2, 3}\sum\limits_{\lambda, \mu=1,2}t_{i\lambda}t_{j\mu}[\theta_i, \, \theta_j]\,\overline{\varphi}_{\lambda}\wedge\overline{\varphi}_{\mu}.\end{equation}

\noindent By (\ref{eqn:hol-vectfiel-comm}), this translates to

\begin{eqnarray}\nonumber\label{eqn:psi_1_brackets2}\frac{1}{2}[\psi_1(t), \, \psi_1(t)] = \frac{1}{2} (t_{11}t_{22}\theta_3\,\overline{\varphi}_1\wedge\overline{\varphi}_2  & + & t_{12}t_{21}\theta_3\,\overline{\varphi}_2\wedge\overline{\varphi}_1 \\
\nonumber & - & t_{21}t_{12}\theta_3\,\overline{\varphi}_1\wedge\overline{\varphi}_2 - t_{22}t_{11}\theta_3\,\overline{\varphi}_2\wedge\overline{\varphi}_1).\end{eqnarray}

\noindent Since $\overline{\varphi}_1\wedge\overline{\varphi}_2 = - \overline{\varphi}_2\wedge\overline{\varphi}_1$, we get

\begin{equation}\label{eqn:psi_1_brackets2}\frac{1}{2}[\psi_1(t), \, \psi_1(t)] = (t_{11}t_{22} - t_{12}t_{21})\,\theta_3\,\overline{\varphi}_1\wedge\overline{\varphi}_2.\end{equation}

\noindent On the other hand, for the choice (\ref{eqn:psi1-def}) we see that

\begin{equation}\label{eqn:d-bar-psi1}\bar\partial\psi_1(t)=d\psi_1(t)=\sum\limits_{i=1}^3\sum\limits_{\lambda=1}^2t_{i\lambda}\,\theta_i\, d\overline{\varphi}_{\lambda}=0\end{equation}

\noindent since $d\overline{\varphi}_1 = d\overline{\varphi}_2=0$. Now setting

\begin{equation}\label{eqn:psi2-def}\psi_2(t):= -(t_{11}t_{22}-t_{12}t_{21})\,\theta_3\overline{\varphi}_3,\end{equation}

\noindent and using (\ref{eqn:Iwa-1forms-rel}) and (\ref{eqn:psi_1_brackets2}), we find

\begin{eqnarray}\label{eqn:d-bar-psi2}\nonumber\bar\partial\psi_2(t) & = & d\psi_2(t) = (t_{11}t_{22}-t_{12}t_{21})\, \theta_3\, (-d\,\overline{\varphi}_3)\\
 & = & (t_{11}t_{22}-t_{12}t_{21})\, \theta_3\,\overline{\varphi}_1\wedge\overline{\varphi}_2 = \frac{1}{2}\, [\psi_1(t), \, \psi_1(t)].\end{eqnarray}

\noindent In particular, $[\psi_1(t), \, \psi_1(t)]$ is seen to be $\bar\partial$-exact here (although it need not be so in the case of an arbitrary manifold, see comments after (\ref{eqn:psi-nu-form})), but the solution $\psi_2(t)$ of equation (\ref{eqn:d-bar-psi2}) need not be of minimal $L^2$-norm (unlike the $\psi_2(t)$ defined in the case of a general manifold by formula (\ref{eqn:psi-nu-form}) for $\nu=2$). In other words, in the special case of the Iwasawa manifold, a solution $\psi_2(t)$ of (\ref{eqn:d-bar-psi2}) is easily observed and we are spared the application of the general formulae (\ref{eqn:psi-nu-form}). This readily yields the desired $\psi(t)$ by setting

\begin{equation}\label{eqn:psi-def}\psi(t):=\psi_1(t) + \psi_2(t) = \sum\limits_{i=1}^3\sum\limits_{\lambda=1}^2t_{i\lambda}\, \theta_i\, \overline{\varphi}_{\lambda} - (t_{11}t_{22} - t_{12}t_{21})\,\theta_3\,\overline{\varphi}_3,\end{equation}

\noindent for which we find 

\begin{equation}\label{eqn:psi-bracket}\frac{1}{2}[\psi(t), \, \psi(t)] = \sum\limits_{j, , k=1}^2\frac{1}{2}[\psi_j(t), \, \psi_k(t)]  = \frac{1}{2}[\psi_1(t), \, \psi_1(t)].\end{equation}

\noindent Indeed, $[\psi_j(t), \, \psi_k(t)] = 0$ for all $(i, \, j)\neq (1, \, 1)$ since these terms involve only brackets of the shape $[\theta_3, \, \theta_i]=0$ and $[\theta_i, \, \theta_3]=0$ which vanish by (\ref{eqn:hol-vectfiel-comm}).

 On the other hand, combining (\ref{eqn:d-bar-psi1}) and (\ref{eqn:d-bar-psi2}), we get

\begin{equation}\label{eqn:d-bar-psi}\bar\partial\psi(t) = \bar\partial\psi_1(t) + \bar\partial\psi_2(t) = \bar\partial\psi_2(t) = \frac{1}{2}[\psi_1(t), \, \psi_1(t)].\end{equation} 

 \noindent Then (\ref{eqn:psi-bracket}) and (\ref{eqn:d-bar-psi}) yield

\begin{equation}\label{eqn:d-bar-intcond}\bar\partial\psi(t) = \frac{1}{2}[\psi(t), \, \psi(t)],\end{equation}

\noindent showing that $\psi(t)$ defined in (\ref{eqn:psi-def}) satisfies the integrability condition (\ref{eqn:int-cond}). 

 By Kuranishi's theorem \ref{The:Kur}, this $T^{1, \, 0}X$-valued $(0, \, 1)$-form $\psi(t)$ defines a {\it locally complete} complex analytic family of deformations $X_t$ of $X$ depending on $6$ {\it effective} parameters $t=(t_{i\lambda})_{1\leq i\leq 3, \, 1\leq\lambda\leq 2}$ such that the complex structure of each fibre $X_t$ is defined by $\bar\partial_t:=\bar\partial - \psi(t)$ and $X_0=X=\C^3/\Gamma$ is the Iwasawa manifold. It is noteworthy that in the special case of the Iwasawa manifold, the power series (\ref{eqn:power-series}) can be built with only two terms ($\psi_1(t)$ and $\psi_2(t)$) and the above simple calculations show $\psi(t)=\psi_1(t)+\psi_2(t)$ to satisfy the integrability condition (\ref{eqn:int-cond}) for all $t=(t_{i\lambda})_{1\leq i\leq 3,\,\,1\leq\lambda\leq 2}\in\Delta_{\varepsilon}\subset\C^6$ if $\varepsilon>0$ is small. With the notation of (\ref{eqn:BvsDelta}), this means that $B=\Delta_{\varepsilon}$.

\vspace{2ex}

  Nakamura goes on to calculate holomorphic coordinates $\zeta_1=\zeta_1(t), \zeta_2=\zeta_2(t), \zeta_3=\zeta_3(t)$ on $X_t$ such that $\zeta_{\nu}(0)=z_{\nu}$ for $\nu=1, 2, 3$ starting from arbitrary holomorphic coordinates $z_1, z_2, z_3$ given beforehand on the Iwasawa manifold $X_0=X=\C^3/\Gamma$. Here is the way he proceeds.

 We are looking for $C^{\infty}$ functions $\zeta_{\nu}(t)$, $\nu=1, 2, 3$, on $X$ satisfying the holomorphicity condition

\begin{equation}\label{eqn:d-bareqn-coord}\bar\partial_t\zeta_{\nu}(t)=0  \,\,\, \Longleftrightarrow \,\,\, \bar\partial\zeta_{\nu}(t) - \psi(t)\zeta_{\nu}(t)=0, \hspace{2ex} \nu=1, 2, 3.\end{equation}

\noindent Given the definition (\ref{eqn:psi-def}) of $\psi(t)$ and the formulae (\ref{eqn:Iwa-hol-vectfiel}) for $\theta_1, \theta_2, \theta_3$, condition (\ref{eqn:d-bareqn-coord}) reads for $\nu=1, 2, 3$\!\!:

\begin{eqnarray}\label{eqn:int-cond-expl}\nonumber\bar\partial\zeta_{\nu} & - & \sum\limits_{\lambda=1}^2t_{1\lambda}\frac{\partial\zeta_{\nu}}{\partial z_1}\, d\bar{z}_{\lambda} - \sum\limits_{\lambda=1}^2t_{2\lambda}\bigg(\frac{\partial\zeta_{\nu}}{\partial z_2} + z_1 \frac{\partial\zeta_{\nu}}{\partial z_3}\bigg) \, d\bar{z}_{\lambda} \\
 & - & \sum\limits_{\lambda=1}^2t_{3\lambda}\frac{\partial\zeta_{\nu}}{\partial z_3}\, d\bar{z}_{\lambda}  + (t_{11}t_{22}-t_{12}t_{21})\frac{\partial\zeta_{\nu}}{\partial z_3} (d\bar{z}_3 - \bar{z}_1 d\bar{z}_2)=0.\end{eqnarray}

 For $\nu=1$, we arrange to have $\frac{\partial\zeta_1}{\partial z_1} = 1$ (in order to get $\zeta_1(t)=z_1 + (\mbox{terms depending only on}\,\, \bar{z}_{\lambda})$) and $\frac{\partial\zeta_1}{\partial z_2} = \frac{\partial\zeta_1}{\partial z_3} =0$. With these choices, condition (\ref{eqn:int-cond-expl}) for $\nu=1$ becomes:

$$\bar\partial\zeta_1 - \sum\limits_{\lambda=1}^2t_{1\lambda} \bar\partial\bar{z}_{\lambda} =0 \hspace{2ex} \Longleftrightarrow \hspace{2ex} \bar\partial\zeta_1(t) = \bar\partial\bigg(\sum\limits_{\lambda=1}^2t_{1\lambda}\bar{z}_{\lambda}\bigg).$$

\noindent Thus we can take

\begin{equation}\label{eqn:zeta1-formula}\zeta_1(t) = z_1 + \sum\limits_{\lambda=1}^2t_{1\lambda}\bar{z}_{\lambda}.\end{equation}

 For $\nu=2$, we similarly require $\frac{\partial\zeta_2}{\partial z_2} = 1$ and $\frac{\partial\zeta_2}{\partial z_1} = \frac{\partial\zeta_2}{\partial z_3} =0$ and condition (\ref{eqn:int-cond-expl}) for $\nu=2$ similarly yields\!\!:

\begin{equation}\label{eqn:zeta2-formula}\zeta_2(t) = z_2 + \sum\limits_{\lambda=1}^2t_{2\lambda}\bar{z}_{\lambda}.\end{equation}

 For $\nu=3$, we require $\frac{\partial\zeta_3}{\partial z_3} = 1$, $\frac{\partial\zeta_3}{\partial z_2} = 0$ and $\frac{\partial\zeta_3}{\partial z_1} = \sum\limits_{\lambda=1}^2 t_{2\lambda}\bar{z}_{\lambda}.$ With these choices, (\ref{eqn:int-cond-expl}) for $\nu=3$ reads

\begin{eqnarray}\label{eqn:int-cond-expl3}\nonumber\bar\partial\zeta_3 & - & \bigg(\sum\limits_{\lambda=1}^2t_{1\lambda}d\bar{z}_{\lambda}\bigg) \bigg(\sum\limits_{\lambda=1}^2t_{2\lambda}\bar{z}_{\lambda}\bigg) - z_1 \sum\limits_{\lambda=1}^2t_{2\lambda}d\bar{z}_{\lambda} \\
\nonumber & - & \sum\limits_{\lambda=1}^2t_{3\lambda}d\bar{z}_{\lambda} + (t_{11}t_{22} - t_{12}t_{21})(d\bar{z}_3 - \bar{z}_1d\bar{z}_2)=0.\end{eqnarray}

\noindent We thus get

\begin{equation}\label{eqn:zeta3-formula}\zeta_3(t) = z_3 + \sum\limits_{\lambda=1}^2(t_{3\lambda} + t_{2\lambda}z_1) \bar{z}_{\lambda} + A(t, \, \bar{z}) - D(t) \bar{z}_3,\end{equation}

\noindent where we have denoted $A(t, \, \bar{z}):= \frac{1}{2}(t_{11}t_{21}\bar{z}_1^2 + 2t_{11}t_{22}\bar{z}_1\bar{z}_2 + t_{12}t_{22}\bar{z}_2^2)$ and $D(t):=(t_{11}t_{22} - t_{12}t_{21})$. We clearly have

$$d\zeta_1(t)\wedge d\zeta_2(t)\wedge d\zeta_3(t) = C(t)\, dz_1\wedge dz_2\wedge dz_3$$

\noindent for a constant $C(t)$ depending in a $C^{\infty}$ way on $t$ such that $C(0)=1$. Hence $\zeta_1(t), \zeta_2(t), \zeta_3(t)$ define holomorphic coordinates on $X_t$ for all $t=(t_{i\lambda})_{1\leq i\leq 3, \, 1\leq\lambda\leq 2}$ such that $\sum\limits_{i=1, 2, 3; \, \lambda=1, 2}|t_{i\lambda}|<\varepsilon$ if $\varepsilon >0$ is small enough.

\vspace{3ex}

\noindent {\bf The example of Alessandrini and Bassanelli proving Theorem \ref{The:bal-notopen}}  \\ 

 In the $6$-parameter Kuranishi family $(X_t)_{t\in B}$, $t=(t_{i\lambda})_{1\leq i\leq 3, \, 1\leq\lambda\leq 2}$, of the Iwasawa manifold $X_0=X=\C^3/\Gamma$, Alessandrini and Bassanelli [AB90] single out the direction corresponding to parameters $t$ such that

\begin{equation}\label{eqn:AB-6thparam}t_{12}\neq 0, \hspace{3ex} t_{ij}=0 \hspace{2ex} \mbox{for all}\,\, (i, \, j)\neq (1, \, 2).\end{equation}

\noindent With this choice of $t$, they have

$$A(t, \, \bar{z})=0  \hspace{2ex}  \mbox{and} \hspace{2ex} D(t)=0.$$

\noindent Thus, denoting $t:=t_{12}$, the holomorphic coordinates of $X_t$ calculated in (\ref{eqn:zeta1-formula}), (\ref{eqn:zeta2-formula}) and (\ref{eqn:zeta3-formula}) reduce to

\begin{equation}\label{eqn:zeta-red-formulae}\zeta_1(t)=z_1 + t\bar{z}_2, \hspace{2ex} \zeta_2(t)=z_2, \hspace{2ex} \zeta_3(t)=z_3.\end{equation} 

\noindent Implicitly $z_1=\zeta_1(t)-t\overline{\zeta}_2(t)$, which yields

\begin{eqnarray}\label{eqn:varphi-new-formulae}\nonumber\varphi_3(t): & = & dz_3 -z_1dz_2 = d\zeta_3(t) + (t\overline{\zeta}_2(t) - \zeta_1(t))\, d\zeta_2(t),\\
 \varphi_2(t): & = & dz_2=d\zeta_2(t), \hspace{3ex} \widetilde{\varphi}_1(t):= dz_1=d\zeta_1(t) - td\overline{\zeta}_2(t).\end{eqnarray}

\noindent Set

\begin{equation}\label{eqn:varphi1-new-formula}\varphi_1(t):=d\zeta_1(t).\end{equation}

\noindent The above $1$-forms $\varphi_1(t), \varphi_2(t), \varphi_3(t)$ are all of $J_t$-type $(1, \, 0)$ since $\zeta_1(t), \zeta_2(t),$ $\zeta_3(t)$ are holomorphic coordinates for the complex structure $J_t$ of $X_t$.

\begin{Prop}(Alessandrini-Bassanelli [AB90, p. 1062]\label{Prop:A-B_examp} Let $(X_t)_t$ be the Kuranishi family of the Iwasawa manifold $X=X_0$, $t=(t_{i\lambda})_{1\leq i\leq 3, \, 1\leq\lambda\leq 2}$. 

 Then, for parameters such that $t_{i\lambda}=0$ for all $(i, \, \lambda)\neq (1, \, 2)$, $X_t$ is {\bf not balanced} for any $t:=t_{12}\neq 0$ satisfying $|t_{12}|<\varepsilon$ if $\varepsilon >0$ is small enough.

\end{Prop}

\noindent {\it Proof.} For the forms defined in (\ref{eqn:varphi-new-formulae}) and (\ref{eqn:varphi1-new-formula}), an immediate calculation shows

\begin{equation}\label{eqn:dphi3-calc} d\varphi_3(t) = (t\,d\bar{\zeta}_2(t) - d\zeta_1(t))\wedge d\zeta_2(t) = -t\,\varphi_2(t)\wedge\overline{\varphi}_2(t) - \varphi_1(t)\wedge\varphi_2(t).\end{equation}

\noindent Thus the $2$-form $d\varphi_3(t)$ has two components: $-t\,\varphi_2(t)\wedge\overline{\varphi}_2(t)$ is of $J_t$-type $(1, \, 1)$, while $- \varphi_1(t)\wedge\varphi_2(t)$ is of $J_t$-type $(2, \, 0)$.

 Recall that $\mbox{dim}_{\C}X_t=3$ for all $t$. Suppose that $X_t$ were balanced for some $t=t_{12}\neq 0$ satisfying $|t_{12}|<\varepsilon$ with $\varepsilon >0$ small. Then there would exist a balanced metric $\omega_t>0$ on $X_t$. Thus $\Omega_t:=\omega_t^2$ would be a $C^{\infty}$ $(2, \, 2)$-form on $X_t$ satisfying

\begin{equation}\label{eqn:bal-square}\Omega_t>0, \hspace{3ex} d\Omega_t=0.\end{equation}

\noindent In this case we would have:

\begin{equation}\label{eqn:zero-neg-contr}0=\int\limits_{X_t}d\Omega_t\wedge i\bar{t}\varphi_3(t) = -\int\limits_{X_t}\Omega_t\wedge i\bar{t}\, d\varphi_3(t) = |t|^2 \int\limits_{X_t}\Omega_t\wedge i\varphi_2(t)\wedge\overline{\varphi_2(t)}.\end{equation}

\noindent Indeed, the first identity above follows from $d\Omega_t=0$ (cf. (\ref
{eqn:bal-square})), the second identity follows from Stokes's theorem, while the third identity follows from formula (\ref{eqn:dphi3-calc}) for $d\varphi_3(t)$ and the fact that the $(2, \, 0)$-component $- \varphi_1(t)\wedge\varphi_2(t)$ is annihilated when wedged with the $(2, \, 2)$-form $\Omega_t$.

 Now $\Omega_t>0$ and $i\varphi_2(t)\wedge\overline{\varphi_2(t)}\geq 0$, hence $\Omega_t\wedge i\varphi_2(t)\wedge\overline{\varphi_2(t)}\geq 0$ at every point of $X_t$. It follows that the right-hand term in (\ref{eqn:zero-neg-contr}) is non-negative. However, since it must vanish by the first identity in (\ref{eqn:zero-neg-contr}), the $(3, \, 3)$-form $\Omega_t\wedge i\varphi_2(t)\wedge\overline{\varphi_2(t)}$ must vanish identically on $X_t$, hence so must the $(1, \, 1)$-form $i\varphi_2(t)\wedge\overline{\varphi_2(t)}$. This can only happen if $\varphi_2(t)$ vanishes identically on $X_t$, which is impossible since $\varphi_2(t)=d\zeta_2(t)$ and $\zeta_2(t)$ is a holomorphic coordinate on $X_t$ if $\varepsilon$ is small enough. This provides the desired contradiction.  

 Therefore $X_t$ cannot be balanced for any $t=t_{12}\neq 0$ if $t_{i\lambda}=0$ for all $(i, \, \lambda)\neq (1, \, 2)$ and $\varepsilon>0$ is small. The proof is complete. \hfill $\Box$

\vspace{2ex}

 It is by means of this Proposition \ref{Prop:A-B_examp} that Alessandrini and Bassanelli proved Theorem \ref{The:bal-notopen}\!\!: they observed that the fibres along one particular direction among the $6$ directions available in the base space of the Kuranishi family of the Iwasawa manifold prove the non-openness of the balanced property under holomorphic deformations.

\vspace{2ex}

 We now make the following

\begin{Obs}\label{Obs:Kurfam-noddbarlemma} (implicit in [Nak75]) In the Kuranishi family of the Iwasawa manifold, the Fr\"olicher spectral sequence does not degenerate at $E_1$ (hence the $\partial\bar\partial$-lemma does not hold) on any fibre $X_t$ corresponding to parameters such that $t_{i\lambda}=0$ for all $(i,\lambda)\neq (1, \, 2)$ and $t:=t_{12}$ satisfies $|t_{12}|<\varepsilon$ with $\varepsilon>0$ small enough.

\end{Obs}

\noindent {\it Proof.} We have seen in (\ref{eqn:H2-Iwa}) that the second Betti number of the Iwasawa manifold is $b_2=8$. By $C^{\infty}$ triviality of the family, all the fibres have the same Betti numbers. On the other hand, Nakamura concludes from his calculations reproduced above (via standard reasoning like the one exemplified above between formulae (\ref{eqn:Iwa-1forms-rel}) and (\ref{eqn:H2-Iwa})) that the Hodge numbers of weight $2$ of any fibre $X_t$ corresponding to parameters such that $t_{i\lambda}=0$ for all $(i,\lambda)\neq (1, \, 2)$ and $t:=t_{12}\neq 0$ are (see [Nak75, table on p. 96 for the case $(ii)$ when $D(t)=0$ and $(t_{11}, \, t_{12}, \, t_{21}, \, t_{22})\neq (0, \, 0, \, 0, \, 0)$])\!\!:

$$h^{2, \, 0}(t)= h^{0, \, 2}(t)=2, \hspace{2ex} h^{1, \, 1}(t)=5.$$

\noindent Thus we see that for any such fibre $X_t$ with $t:=t_{12}\neq 0$, we have\!\!:

$$b_2=8< h^{2, \, 0}(t) + h^{1, \, 1}(t) + h^{0, \, 2}(t) =9.$$

 Hence the conclusion follows. Notice that for fibres as above with $t:=t_{12}\neq 0$, Nakamura's table gives $b_1=4 = 2 + 2 = h^{1, \, 0}(t) + h^{0, \, 1 }(t)$, while we have seen in Observation \ref{Obs:coh-Iwasawa} that for $X_0$ we have $b_1=4 < 3 + 2 = h^{1, \, 0}(0) + h^{0, \, 1 }(0)$.  \hfill $\Box$

\vspace{3ex}

\noindent {\bf The desired examples of sG manifolds}\\

 We can now conclude this section by exhibiting the desired examples of sG manifolds showing the difference between, on the one hand, the sG property and, on the other hand, the (combined) balanced and $\partial\bar\partial$-lemma properties. It suffices to bring together Theorem \ref{The:sG-openness}, Theorem \ref{The:bal-notopen} and Observation \ref{Obs:Kurfam-noddbarlemma}.

\begin{The}\label{The:sG-examp} Let $(X_t)_t$ be the Kuranishi family of the Iwasawa manifold $X=X_0$, $t=(t_{i\lambda})_{1\leq i\leq 3, \, 1\leq\lambda\leq 2}$. 

 Then, for parameters such that $t_{i\lambda}=0$ for all $(i, \, \lambda)\neq (1, \, 2)$, $X_t$ is a {\bf strongly Gauduchon} manifold that is {\bf not balanced} and whose {\bf Fr\"olicher spectral sequence does not degenerate at $E_1$} (hence the {\bf $\partial\bar\partial$-lemma does not hold}) for any $t=t_{12}\neq 0$ satisfying $|t_{12}|<\varepsilon$ if $\varepsilon >0$ is small enough.

\end{The}

\noindent {\it Proof.} Since the Iwasawa manifold is balanced (cf. Corollary \ref{Cor:2}), it is also an sG manifold. Since the sG property is open under holomorphic deformations (cf. Theorem \ref{The:sG-openness}), all sufficiently nearby fibres $X_t$ in the Kuranishi family of the Iwasawa manifold $X_0$ are again sG manifolds. However, by the observation of Alessandrini and Bassanelli (cf. Proposition \ref{Prop:A-B_examp}), the fibres $X_t$ corresponding to parameters for which $t_{i\lambda}=0$ for all $(i, \, \lambda)\neq (1, \, 2)$ are not balanced if $t:=t_{12}$ is sufficiently close to $0$. By Observation \ref{Obs:Kurfam-noddbarlemma}, the Fr\"olicher spectral sequence does not degenerate at $E_1$, hence the $\partial\bar\partial$-lemma does not hold, on any of these fibres. \hfill $\Box$

\section{The Eastwood-Singer construction}\label{section:ES}

 In this section we give an outline of the Eastwood-Singer proof [ES93] of part $(b)$ of Theorem \ref{The:E1=Einf_op-cl} asserting the non-closedness of the Fr\"olicher degeneration at $E_1$ under holomorphic deformations. All the fibres of the holomorphic family they construct in their example are twistor spaces (hence in particular compact complex manifolds of dimension $3$). 

 We briefly recall the barest essentials of Penrose's {\it twistor space} theory [Pen76] for which the standard mathematical reference is [AHS78]. Given a compact oriented Riemannian manifold $(M, \, g)$ of real dimension $4$, for every point $x\in M$ one defines $L_x$ to be the set of all complex structures $J_x$ on the tangent space $T_xM$ such that $J_x$ is orthogonal w.r.t. $g_x$ and $J_x$ induces the negative (i.e. opposite to the given one) orientation on $T_xM$. Orthogonality w.r.t. $g_x$ for a complex structure $J_x : T_xM\to T_xM$ means, as usual, that

$$g_x(J_xu, \, J_xv)=g_x(u, \, v), \hspace{2ex} \mbox{for all} \,\, u, v\in T_xM,$$

\noindent i.e. $J_x$ is required to be a $g_x$-isometry of $T_xM$. When $x$ varies in $M$, the union $Z$ (which depends only on the conformal class $[g]$ of Riemannian metrics on $M$ but not on the actual representative $g$) of all twistor lines $L_x$ has a natural structure as a $C^{\infty}$ manifold of real dimension $6$ and a natural almost complex structure. The almost complex structure is {\it integrable} if and only if the conformal structure $[g]$ of $M$ is {\it self-dual}. Recall that for an arbitrary-dimensional Riemannian manifold $(M, \, g)$, the curvature tensor $R$ decomposes as

$$R=W + \rho,$$

\noindent where $W$ is the Weyl tensor depending only on the conformal class $[g]$. Peculiar to the case when $M$ has real dimension $4$ is a further decomposition of the Weyl tensor as

$$W=W_{+} + W_{-},$$

\noindent where $W_{+}$ is the {\it self-dual} component and $W_{-}$ is the {\it anti-self-dual} component. Reversing the orientation of $M$ permutes $W_{+}$ and $W_{-}$. The conformal structure $[g]$ of $M$ is said to be {\it self-dual} if the associated Weyl tensor reduces to its {\it self-dual} component $W_{+}$ (i.e. $W_{-}=0$).

 One thus gets the Penrose correspondence between {\it self-dual} compact connected oriented $C^{\infty}$ real $4$-manifolds $M$ and the associated twistor spaces $Z$ (which are compact complex $3$-manifolds). The natural projection $\nu : Z\to M$ can be identified with the unit sphere subbundle of the rank-three real vector bundle of anti-self-dual $2$-forms on $M$. Every twistor line $\nu^{-1}(x)=L_x$ is isomorphic to the complex projective line $\Proj^1$.  

 For every $r$, denote by $\Lambda^r:=C^{\infty}_r(M, \, \C)$ the space of $C^{\infty}$ complex-valued $r$-forms on $M$. Recall that the Hodge star operator of the Riemannian metric $g$ of $M$ acting on $2$-forms

$$\star : \Lambda^2\to\Lambda^2$$

\noindent satisfies $\star^2=1$. Hence it induces a direct-sum splitting

$$\Lambda^2=\Lambda^2_{+} \oplus \Lambda^2_{-1}$$

\noindent into its $\pm 1$-eigenspaces. The $2$-forms $u\in\Lambda^2_{+}$ (i.e. $\star\, u=u$) are termed {\it self-dual}, while the $2$-forms $u\in\Lambda^2_{-}$ (i.e. $\star\, u=-u$) are termed {\it anti-self-dual}. If the Weyl curvature tensor $W$ is viewed as a bundle-valued $2$-form, its components satisfy $\star\, W_{+}=W_{+}$ and respectively $\star\, W_{-}=-W_{-}$. One gets a corresponding splitting of the differential operator $d$ acting on $1$-forms of $M$\!\!:

$$\Lambda^1\xrightarrow{d=d_{+} + d_{-}}\Lambda^2=\Lambda^2_{+} \oplus \Lambda^2_{-}$$

\noindent that induces cohomology groups $H^2_{+}(M, \, \C)$ and $H^2_{-}(M, \, \C)$.

\vspace{2ex}

 For a fixed compact connected oriented $C^{\infty}$ manifold $M$ of real dimension $4$ endowed with a {\it self-dual} Riemannian metric $g$, Eastwood and Singer establish the following general facts about the associated twistor space $Z$. The details are to be found in [ES93, $\S.\,\, 2,3,4.$]. \\

\noindent {\bf Fact 1.} Given that the sheaves $\Omega^1$, $\Omega^2$ and $\Omega^3$ of germs of holomorphic $1$, $2$ and respectively $3$-forms on $Z$ are explicitly given, on every twistor line $L=L_x=\Proj^1$, by the formulae

$$\Omega^1_{|L}\simeq{\cal O}(-2)\oplus{\cal O}(-1)\oplus{\cal O}(-1), \hspace{2ex} \Omega^2_{|L}\simeq{\cal O}(-3)\oplus{\cal O}(-3)\oplus{\cal O}(-2), \hspace{2ex} \Omega^3_{|L}\simeq{\cal O}(-4),$$

\noindent we see that the restricted bundles $\Omega^1_{|L}, \Omega^2_{|L}, \Omega^3_{|L}$ have no non-trivial sections over $L$. Since $Z$ is fibred by projectives lines $L=L_x=\Proj^1$, it follows that the vector bundles $\Omega^1, \Omega^2, \Omega^3$ have no non-trivial global holomorphic sections over $Z$, hence

\begin{equation}\label{eqn:D-coh-nosect}H^{1, \, 0}(Z, \, \C)=0, \hspace{2ex} H^{2, \, 0}(Z, \, \C)=0, \hspace{2ex} H^{3, \, 0}(Z, \, \C)=0.\end{equation}

\noindent By Serre duality, one also gets

\begin{equation}\label{eqn:D-coh-Sduality}H^{2, \, 3}(Z, \, \C)=0, \hspace{2ex} H^{1, \, 3}(Z, \, \C)=0, \hspace{2ex} H^{0, \, 3}(Z, \, \C)=0.\end{equation}

\noindent Since $Z$ and $M$ are compact, one infers that

\begin{equation}\label{eqn:D-coh-comp}H^{0, \, 0}(Z, \, \C)= H^0(M, \, \C)=\C \hspace{2ex} \mbox{and} \hspace{2ex} H^{3, \, 3}(Z, \, \C)= H^4(M, \, \C)=\C,\end{equation}

\noindent where the latter set of identities follows from the former by Serre duality on $Z$ and Poincar\'e duality on $M$.

\vspace{2ex}

\noindent {\bf Fact 2.} The Penrose transform relating analytic cohomology on the twistor space $Z$ to solutions of differential equations on the base manifold $M$ is used to see that \\

\noindent $\bullet$ for all $q=0, 1, 2, 3$, $H^{0, \, q}(Z, \, \C)$ is canonically isomorphic to the $q^{th}$ cohomology of the complex 

$$\Lambda^0\longrightarrow\Lambda^1\xrightarrow{d_{-}}\Lambda^2_{-}\longrightarrow 0.$$

\noindent $\bullet$ for all $q=0, 1, 2, 3$, $H^{3, \, q}(Z, \, \C)$ is canonically isomorphic to the $q^{th}$ cohomology of the complex 

$$0\longrightarrow\Lambda^2_{-}\stackrel{d}{\longrightarrow}\Lambda^3\longrightarrow\Lambda^4.$$

 It follows as a corollary that

\begin{eqnarray}\label{eqn:fact2-coh}\nonumber H^{0, \, 1}(Z, \, \C)=H^1(M, \, \C) & , & H^{0, \, 2}(Z, \, \C)=H^2_{-}(M, \, \C)\\
H^{3, \, 1}(Z, \C)=H^2_{-}(M, \, \C) & , & H^{3, \, 2}(Z, \, \C)=H^3(M, \, \C).\end{eqnarray}

\noindent {\bf Fact 3.} The Penrose transform of the vector bundles associated with $\Omega^1$ and $\Omega^2$ yields a commutative diagram with exact rows (cf. [ES93, Proposition 3.3]) from which it follows that\!\!: 



\vspace{2ex}

\noindent $\bullet$ $H^{1, \, 1}(Z, \, \C)$ can be identified with the set of pairs $(f, \, \rho)\in\Lambda^0\oplus\Lambda^2_{+}$ satisfying the equation

\begin{equation}\label{eqn:H11-eqn}{\cal D}df - d\rho=0  \hspace{2ex}\mbox{on}\,\, M,\end{equation}

\noindent where the operator ${\cal D} : \Lambda^1\longrightarrow\Lambda^3$ is defined by

\begin{equation}\label{eqn:def-calD}\omega_b\mapsto \bigg(\nabla^a\nabla^b + 2R^{ab}-\frac{2}{3}R\, g^{ab}\bigg)\,\omega_b,\end{equation}

\noindent while $\nabla$ denotes the Levi-Civita connection, $R_{ab}$ denotes the Ricci curvature and $R$ denotes the scalar curvature of $(M, \, g)$; 

\vspace{1ex}

\noindent $\bullet$ For $E^{1, \, 1}_2$ featuring in the Fr\"olicher spectral sequence of $Z$ at $E_2$ level, we always have the canonical isomorphism of $\C$-vector spaces: 

\begin{equation}\label{eqn:E11_2-ident}E^{1, \, 1}_2:=\ker\bigg(d_1=\partial : H^{1, \, 1}(Z, \, \C)\rightarrow H^{2, \, 1}(Z, \, \C)\bigg)\simeq H^0(M, \, \C)\oplus H^2_{+}(M, \C),\end{equation}

\noindent hence $\mbox{dim}_{\C}E^{1, \, 1}_2 = 1 + b_{+}(M),$ where $b_{+}(M)\!\!:=\mbox{dim}_{\C}H^2_{+}(M, \, \C)$. (Clearly $\mbox{dim}_{\C}H^0(M, \, \C)=1$ by compactness of $M$.)

\vspace{2ex}

\noindent {\bf Fact 4.} By [Hit81], there exists a $(1, \, 1)$-form $h$ on $Z$ such that $dh=0$ and $h_{|L}$ generates $H^2(L, \, \C)$ for every twistor line $L=L_x=\nu^{-1}(x)\subset Z$. Using the Leray-Hirsch theorem, we get isomorphisms

$$H^{r-2}(M, \, \C)\oplus H^r(M, \, \C)\simeq H^r(Z, \, \C), \hspace{3ex} 0\leq r\leq 6,$$

\noindent induced by $\Lambda^{r-2}\oplus\Lambda^r\ni (\alpha, \, \beta)\mapsto h\wedge\nu^{\star}\alpha + \nu^{\star}\beta\in C^{\infty}_r(Z, \, \C)$.

\vspace{2ex}

 The conclusion of these four facts is that for any twistor space $Z$ we always have $E_2(Z)=E_{\infty}(Z)$ (i.e. the Fr\"olicher spectral sequence degenerates {\it at the latest} at $E_2$ level), while the degeneration at $E_1$ level depends exclusively on one arrow. Recall that for any complex manifold, the $E_1$ level of the Fr\"olicher spectral sequence is given by the Dolbeault cohomology groups (i.e. $E_1^{p, \, q}=H^{p, \, q}$ for all $p, q$) and the $d_1$ arrows are induced by $\partial$:

$$E_1^{p, \, q}\xrightarrow{d_1=\partial}E_1^{p+1, \, q} \hspace{2ex}\mbox{for all}\,\, p, q.$$

\noindent When $Z$ is a twistor space, it follows from the above facts $(1)-(4)$ that all arrows $d_1=\partial : E_1^{p, \, 0}\rightarrow E_1^{p+1, \, 0}$ are zero when $q=0$ (hence also when $q=3$ by duality), while the part of the $E_1$ level of the Fr\"olicher spectral sequence corresponding to $q=1$ reads\!\!: \\

\noindent $\begin{array}{lllll}\label{eqn:Fsq-1} 0\rightarrow & H^{0, \, 1}(Z, \, \C) & \xrightarrow{d_1=0}H^{1, \, 1}(Z, \, \C)\stackrel{d_1}{\rightarrow}H^{2, \, 1}(Z, \, \C)\xrightarrow{d_1=0} & H^{3, \, 1}(Z, \, \C)  & \rightarrow 0\\
 & \parallel &  & \parallel &  \\
 & H^1(M, \, \C) &  & H^2_{-}(M, \, \C)      &   .\end{array}$ \\

\noindent Since the part corresponding to $q=2$ at $E_1$ level is dual to that for $q=1$, we see that $E_1(Z)=E_{\infty}(Z)$ if and only if the middle arrow above, i.e.

\begin{equation}\label{eqn:midarrow}d_1=\partial : H^{1, \, 1}(Z, \, \C)\to H^{2, \, 1}(Z, \, \C),\end{equation}

\noindent vanishes. It is clear that the kernel of the arrow (\ref{eqn:midarrow}) equals $E_2^{1, \, 1}$ (cf. (\ref{eqn:E11_2-ident})) since the arrow preceding it in the last complex vanishes. We get the

\begin{Conc}(cf. [ES93, p. 653-662])\label{Conc:E1-deg} For any twistor space $\nu : Z\to M$, $H^0(M, \, \C)\oplus H^2_{+}(M, \, \C)$ injects canonically into $H^{1, \, 1}(Z, \, \C)$. 

 Moreover, the Fr\"olicher spectral sequence of $Z$ degenerates at $E_1$ if and only if $H^{1, \, 1}(Z, \, \C)$ is isomorphic to $H^0(M, \, \C)\oplus H^2_{+}(M, \, \C)$ iff $h^{1, \, 1}(Z)=1 + b_{+}(M)$.

\end{Conc}

 To construct the actual examples that prove part $(b)$ of Theorem \ref{The:E1=Einf_op-cl}, Eastwood and Singer [ES93, $\S.\,\, 5$] go on to choose $M$ to be a compact complex surface endowed with a K\"ahler metric $\omega$ of zero scalar curvature. By [Leb86], any such $M$ that has been given the conjugate orientation is self-dual, hence $M$ possesses a twistor space $Z$. 

 On the other hand, for any compact K\"ahler manifold $(X, \, \omega)$, the {\it Lichnerowicz operator} (cf. e.g. [Bes87]) is defined on functions by

\begin{equation}\label{eqn:L-op}{\cal L}:C^{\infty}(X, \, \C)\to C^{\infty}(X, \, \C), \hspace{2ex} {\cal L}(f):=\Delta^2f + \langle\langle dd^cf, \, \mbox{Ric}\,\omega\rangle\rangle,\end{equation}

\noindent where $\Delta$ is the Laplacian and $\mbox{Ric}\,\omega$ is the Ricci form of $\omega$. A result of Lichnerowicz (cf. e.g. [Bes87, Proposition 2.151]) guarantees that when the scalar curvature of $(X, \, \omega)$ is constant, there is an isomorphism

$$u(X)/u_0(X)\simeq (\ker{\cal L})/\C,$$

\noindent where $u(X)$ denotes the complex Lie algebra of holomorphic vector fields on $X$, while $u_0(X)$ denotes the abelian Lie algebra of parallel such vector fields.

 It turns out that when $(X, \, \omega)=(M, \, \omega)$ is a compact K\"ahler complex surface of zero scalar curvature, the operator ${\cal D}$ of (\ref{eqn:def-calD}) (defining equation (\ref{eqn:H11-eqn}) on $M$ which characterises the Dolbeault cohomology group $H^{1, \, 1}(Z, \, \C)$) relates to the Lichnerowicz operator as follows (cf. [ES93, p. 662])\!\!:

$$d{\cal D}d={\cal L},$$

\noindent after identification of $4$-forms with functions on $M$ via the volume form.

 This leads to the following consequence of the main theorem of [ES93].

\begin{The}\label{The:ES-mainth-conseq}([ES93, Theore 5.3.]) Let $(M, \, \omega)$ be a compact K\"ahler complex surface of zero scalar curvature. If $Z$ is the twistor space of $M$, then

\begin{equation}\label{eqn:H11-finform}\nonumber H^{1, \, 1}(Z, \, \C)=H^0(M, \, \C)\oplus H^2_{+}(M, \,\C)\oplus\frac{\{\mbox{holomorphic vector fields on}\,\, M\}}{\{\mbox{parallel vector fields on}\,\, M\}}.\end{equation}

 Since we always have $E_2^{1, \, 1}\simeq H^0(M, \, \C)\oplus H^2_{+}(M, \,\C)=\C\oplus H^2_{+}(M, \,\C)$ (cf. (\ref{eqn:E11_2-ident})), we see by Conclusion \ref{Conc:E1-deg} that the Fr\"olicher spectral sequence of $Z$ degenerates at $E_1$ if and only if all holomorphic vector fields on $M$ are parallel.

\end{The}

 The simplest case to which Eastwood and Singer apply their results is that of $M:=\Sigma_g\times\Proj^1$, where $\Sigma_g$ is any compact complex curve of genus $g\geq 2$ endowed with the Poincar\'e metric $\omega_P$ (which is of constant curvature $-1$) and the complex projective line $\Proj^1$ is endowed with the metric $\omega_{\Proj^1}:=2\,\omega_{FS}$ (which is of constant curvature $+1$ if $\omega_{FS}$ denotes the Fubini-Study metric). The metric $\omega$ induced on the product compact complex surface $M=\Sigma_g\times\Proj^1$ is a K\"ahler metric of zero scalar curvature.  

\begin{Prop}\label{Prop:Frolicher-nondeg}([ES93, p. 663]) If $M=\Sigma_g\times\Proj^1$ has been given the K\"ahler metric of zero scalar curvature $\omega=\omega_P\oplus\omega_{\Proj^1}$, the Fr\"olicher spectral sequence of the twistor space $Z$ of $M$ does not degenerate at $E_1$.

\end{Prop}

\noindent {\it Proof.} We have seen in (\ref{eqn:H11-eqn}) that 

\begin{equation}\label{eqn:re-H11-eqn}E_1^{1, \, 1}=H^{1, \, 1}(Z, \, \C)\simeq\{(f, \, \rho)\in\Lambda^0\oplus\Lambda_{+}^2\,\,;\,\,{\cal D}df=d\rho\ \,\,\mbox{on}\,\, M\}.\end{equation} 

\noindent On $M=\Sigma_g\times\Proj^1$, we have

$${\cal D}df = (\nabla^a\nabla^b + 2R^{ab})\nabla_bf$$

\noindent because the scalar curvature $R\equiv 0$. If a $C^{\infty}$ function $f: \Sigma_g\times\Proj^1\to \C$ depends only on the $\Proj^1$ variable, we get:

\begin{eqnarray}\label{eqn:2-laplacian}\nonumber{\cal D}df & = & (\nabla^i\nabla^j + 2\omega^{ij})\nabla_jf = \nabla^i(\nabla^j\nabla_jf) + 2\omega^{ij}\nabla_jf\\
\nonumber & = & \nabla^i(-\Delta + 2)f,\end{eqnarray}

\noindent where $\nabla_j$ is the Chern connection and $\Delta=-\nabla^j\nabla_j$ is the Laplacian on $\Proj^1$, while $\omega=(\omega_{ij})_{i,\, j}$ locally on $\Proj^1$ and the identity $R^{ij}=\omega^{ij}$ holds because $\mbox{Ric}\,\omega=\omega$ on $\Proj^1$. Using now the well-known fact that $2$ is an eigenvalue of the Laplacian on $\Proj^1$ with eigenspace $E_{\Delta}(2)$ of complex dimension $3$, we get 

$$\{C^{\infty} \,\,\mbox{functions}\,\, f : M\rightarrow \C \,\, ;\,\, {\cal D}df=0\}\supset \{\mbox{Constants}\}\oplus E_{\Delta}(2)= \C \oplus E_{\Delta}(2),$$ 

\noindent hence $\mbox{dim}_{\C}\{C^{\infty} \,\,\mbox{functions}\,\, f : M\rightarrow \C \,\, ;\,\, {\cal D}df=0\}\geq 4$. Thus (\ref{eqn:re-H11-eqn}) yields

$$h^{1, \, 1}(Z)=\mbox{dim}_{\C}H^{1, \, 1}(Z, \, \C)\geq 4 + \mbox{dim}_{\C}\{\rho\in\Lambda^2_{+}\,\, ; \,\, d\rho=0\} = 4 + b_{+}(M).$$

\noindent As $\mbox{dim}_{\C}E_2^{1, \, 1} = 1+b_{+}(M)$ by (\ref{eqn:E11_2-ident}), the claim follows from Conclusion \ref{Conc:E1-deg}.  \hfill $\Box$

 \vspace{2ex}

 It is obvious that $M=\Sigma_g\times\Proj^1$ is the ruled surface $\Proj(E_0)$ associated with the trivial rank-two vector bundle $E_0:=\Sigma_g\times\C^2\to\Sigma_g$. This suggests a natural way of constructing a holomorphic family of compact complex surfaces $(M_t)_{t\in\Delta}$ such that $M_0=M=\Sigma_g\times\Proj^1$ for which we can hope to ensure that $E_1(Z_t)=E_{\infty}(Z_t)$ for all $t\in\Delta^{\star}$ in the associated family of twistor spaces $(Z_t)_{t\in\Delta}$: take $M_t:=\Proj(E_t)$ when we have found a {\it suitable} family of rank-two holomorphic vector bundles $(E_t)_{t\in\Delta}$ over $\Sigma_g$ with the trivial bundle $E_0:=\Sigma_g\times\C^2$ corresponding to $t=0$. The authors make clear the meaning of {\it suitable} in the following form.

\begin{Prop}\label{Prop:ES-NSR}([ES93, p. 663-664]) If $E\rightarrow\Sigma_g$ is a {\bf stable} rank-two holomorphic vector bundle with trivial determinant over a compact complex curve of genus $g\geq 2$, then the corresponding ruled surface $\Proj(E)$ satisfies\!:

\noindent $(i)$\, $\Proj(E)$ admits a K\"ahler metric of zero scalar curvature\!;

\noindent $(ii)$\, $\Proj(E)$ has no non-zero holomorphic vector fields.

\end{Prop}

\noindent {\it Proof.} Eastwood and Singer deduce the above statements from classical results of Narasimhan-Seshadri [NS65] and Narasimhan-Ramanan [NR69]. 

 By [NS65], any stable rank-two holomorphic vector bundle on $\Sigma_g (g\geq 2)$ arises from a representation of $\pi_1(\Sigma_g)$ into $SU_2$, hence the corresponding ruled surface is a quotient $\Proj(E)=(\Proj^1\times H)/\pi_1$ where $H$ is the upper half-plane and $\pi_1$ acts by isometries of the natural metric on $\Proj^1\times H$. This natural metric, obtained as the Riemannian product of the metric of curvatute $+1$ on $\Proj^1$ with that of curvature $-1$ on $H$, is K\"ahler and of zero scalar curvature, hence $(i)$ follows.

 It follows from [NS65] that any vector bundle as in the statement satisfies $H^0(\Sigma_g, \, S^2E^{\star})=0$. Indeed, for any rank-two vector bundle $E$ there is an isomorphism of bundles (see e.g. the proof of Proposition 3.3. in [Hit87])\!\!:

\begin{equation}\label{eqn:iso-bundlesS2}\mbox{End}_0E\simeq S^2E^{\star}\otimes\Lambda^2E = S^2E^{\star}\otimes\det E\end{equation}

\noindent induced by the map $T$ which associates with every $A\in\mbox{End}\,E$ the quadratic map $E\ni v\mapsto Av\wedge v\in\Lambda^2E$. The kernel of $T$ consists of the scalar endomorphisms (identified with $\C$), while $\mbox{End}_0E$ denotes the traceless endomorphisms. It obviously satisfies $\mbox{End}\,E = \mbox{End}_0E\oplus\C.$ On the other hand, if $E$ is stable then $E$ is {\it simple} (cf. [NS65, Corollary to Proposition 4.3.]), i.e. $\mbox{dim}_{\C}H^0(\Sigma_g, \, \mbox{End}\,E)=1$, which means that the only endomorphisms of $E$ are the scalar ones. Hence $H^0(\Sigma_g, \, \mbox{End}_0\,E)=0$ when $E$ is stable. Thus, if $E$ is stable and $\det E$ is trivial, we see by (\ref{eqn:iso-bundlesS2}) that $H^0(\Sigma_g, \, S^2E^{\star})=0$.

Now, an easy argument explained in [ES93] shows that $H^0(\Proj(E), \, T^{1, \, 0}\Proj(E))\simeq H^0(\Sigma_g, \, S^2E^{\star})$ for any rank-two holomorphic vector bundle $E\rightarrow\Sigma_g$ when $g\geq 2$. This is because any holomorphic vector field $\xi\in H^0(\Proj(E), \, T^{1, \, 0}\Proj(E))$ must be {\it vertical}. Indeed, the normal bundle of any fibre of $\Proj(E)\rightarrow\Sigma_g$ being trivial, the component of $\xi$ normal to any such fibre is constant along the fibre, hence it projects to a holomorphic vector field on $\Sigma_g$. Now any holomorphic vector field on $\Sigma_g$ must vanish because, since $g\geq 2$, $\Sigma_g$ embeds into a $g$-dimensional complex torus $\C^g/\Lambda$ whose flat metric induces a metric of negative curvature on $T^{1, \, 0}\Sigma_g$ viewed as a holomorphic line subbundle of the tangent bundle of $\C^g/\Lambda$. Thus any $\xi\in H^0(\Proj(E), \, T^{1, \, 0}\Proj(E))$ is indeed vertical, hence the restriction of $\xi$ to any fibre of $\Proj(E)\rightarrow\Sigma_g$ defines a holomorphic vector field on the fibre $\Proj(E_t)\simeq\Proj^1$. Meanwhile, the holomorphic vector fields of $\Proj^1$ are the holomorphic sections of $-K_{\Proj^1}={\cal O}_{\Proj^1}(2)$. Thus $(ii)$ follows.  \hfill $\Box$

\vspace{2ex}

 We now get an immediate corollary of Theorem \ref{The:ES-mainth-conseq} and Proposition \ref{Prop:ES-NSR}. In view of Theorem \ref{The:ES-mainth-conseq}, a weaker version of conclusion $(ii)$ of Proposition \ref{Prop:ES-NSR} (with {\it non-parallel} in place of {\it non-zero}) would have sufficed.

\begin{Cor}\label{Cor:ES-surf-twis-fam}([ES93, p. 664]) Let $E\rightarrow\Sigma_g$ be a {\bf stable} rank-two holomorphic vector bundle with trivial determinant over a compact complex curve of genus $g\geq 2$. Let $Z$ be the twistor space associated with the ruled surface $\Proj(E)$.

 Then the Fr\"olicher spectral sequence of $Z$ degenerates at $E_1$. 

\end{Cor}

 Putting together Proposition \ref{Prop:Frolicher-nondeg} and Corollary \ref{Cor:ES-surf-twis-fam}, we see that it suffices to show that the trivial rank-two vector bundle $E_0=\Sigma_g\times\C^2\rightarrow\Sigma_g$ deforms to {\bf stable} rank-two holomorphic vector bundles $E_t\rightarrow\Sigma_g$ with trivial determinant, $t\in\Delta^{\star}$. As explained before the statement of Proposition \ref{Prop:ES-NSR}, associated with the family of bundles $(E_t)_{t\in\Delta}$ will be the family of ruled surfaces $(M_t:=\Proj(E_t))_{t\in\Delta}$ whose corresponding family of twistor spaces $(Z_t)_{t\in\Delta}$ will provide the example proving the Eastwood-Singer part $(b)$ of Theorem \ref{The:E1=Einf_op-cl}: $E_1(Z_t)=E_{\infty}(Z_t)$ for all $t\neq 0$ by Corollary \ref{Cor:ES-surf-twis-fam}, but $E_1(Z_0) \neq E_{\infty}(Z_0)$ by Proposition \ref{Prop:Frolicher-nondeg}.

\begin{The}\label{The:ES_triv-stab-def}([ES93, p. 664-665]) Given any compact complex curve $\Sigma_g$ of genus $g\geq 2$, the trivial rank-two vector bundle $E_0=\Sigma_g\times\C^2\rightarrow\Sigma_g$ deforms to {\bf stable} rank-two holomorphic vector bundles with trivial determinant on $\Sigma_g$.

\end{The}

\noindent {\it Proof.} We give an outline of the proof found in [ES93] which, as mentioned there, is a modification of an argument from [NR69]. Fix an arbitrary point $x\in\Sigma_g$ and denote by $L_x$ the holomorphic line bundle on $\Sigma_g$ defined by $x$ (viewed as a divisor on the curve $\Sigma_g$). If $L_x^{-1}$
denotes the line bundle dual to $L_x$, consider the rank-two holomorphic vector bundles $E\rightarrow\Sigma_g$ that are non-trivial extensions of $L_x$ by $L_x^{-1}$. Equivalently, one considers short exact sequences of vector bundles on $\Sigma_g$\!:

\begin{equation}\label{eqn:ext-Lx-Lx-1}0\rightarrow L_x^{-1} \rightarrow E \rightarrow L_x \rightarrow 0\end{equation}

\noindent that do not split holomorphically. There is a one-to-one correspondence between the equivalence classes of such non-trivial extensions and the non-zero classes $\{\beta^{\star}\}\in H^{0, \, 1}(\Sigma_g, \, \mbox{Hom}(L_x, \, L_x^{-1}))\simeq H^1(\Sigma_g, \, {\cal O}(L_x^{-2}))$, where $\beta^{\star}$ denotes the second fundamental form of extension (\ref{eqn:ext-Lx-Lx-1}). Any such $E$ has trivial determinant (since $\det E = L_x^{-1}\otimes L_x$) and degree zero since $L_x$ has degree $+1$, $L_x^{-1}$ has degree $-1$ and the degree is additive in exact sequences. By [NR69, Lemma 5.1], any such $E$ can have no holomorphic line subbundles of positive degree, hence any such $E$ is at least {\it semi-stable}. To find {\it stable} vector bundles $E$ arising as extensions (\ref{eqn:ext-Lx-Lx-1}), it remains to rule out the existence of holomorphic line subbundles of degree zero in $E$.

 Suppose that a rank-two holomorphic vector bundle $E$ given by an extension (\ref{eqn:ext-Lx-Lx-1}) contains a degree-zero holomorphic line subbundle $L\subset E$. By the proof of [NR69, Lemma 5.2], any such $L$ must be of the form $L=L_x\otimes L_y^{-1}$ for some point $y\in\Sigma_g$ and Lemma 5.2 in [NR69] asserts that the class $\{\beta^{\star}\}$ of extension (\ref{eqn:ext-Lx-Lx-1}) maps to the zero class in the cohomology group $H^{0, \, 1}(\Sigma_g, \, \mbox{Hom}(L, \, L_x^{-1})) = H^1(\Sigma_g, \, {\cal O}(L_x^{-2}\otimes L_y))$ under the map

\begin{equation}\label{eqn:surj-map}H^1(\Sigma_g, \, {\cal O}(L_x^{-2}))\rightarrow H^1(\Sigma_g, \, {\cal O}(L_x^{-2}\otimes L_y)).\end{equation}

\noindent Since (\ref{eqn:surj-map}) is a surjective linear mapping of complex vector spaces of dimensions $g+1$ and respectively $g$ thanks to Riemann-Roch, the kernel of (\ref{eqn:surj-map}) is a complex line $l_y$ in the $(g+1)$-dimensional vector space $H^1(\Sigma_g, \, {\cal O}(L_x^{-2}))$. We conclude that $L=L_x\otimes L_y^{-1}$ is a holomorphic subbundle of $E$ iff

\begin{equation}\label{eqn:ext-class-kercond}\nonumber\{\beta^{\star}\}\in\ker\bigg(H^1(\Sigma_g, \, {\cal O}(L_x^{-2}))\rightarrow H^1(\Sigma_g, \, {\cal O}(L_x^{-2}\otimes L_y))\bigg)\setminus\{0\}=l_y\setminus\{0\}.\end{equation}

\noindent By Lemma 3.3 in [NR69], proportional extension classes $\{\beta^{\star}\}$ and $\lambda\{\beta^{\star}\}$, with $\lambda\in\C^{\star}$, give rise to isomorphic bundles $E$ and $E_{\lambda}$. Now the punctured line $l_y\setminus\{0\}$ defines a point in the complex projective space $\Proj H^1(\Sigma_g, \, {\cal O}(L_x^{-2}))\simeq\Proj^g$ and, when $y$ varies in $\Sigma_g$, we get an analytic mapping

\begin{equation}\label{eqn:ext-class-projmap}\Sigma_g\ni y\mapsto l_y\setminus\{0\}\in\Proj H^1(\Sigma_g, \, {\cal O}(L_x^{-2}))\simeq\Proj^g.\end{equation}

 The conclusion is that an extension class $\{\beta^{\star}\}\in H^1(\Sigma_g, \, {\cal O}(L_x^{-2}))\setminus\{0\}$ defines a holomorphic rank-two vector bundle $E\rightarrow\Sigma_g$ with trivial determinant that contains no degree-zero holomorphic line subbundles $L\subset E$ (hence $E$ is {\it stable}) if and only if the image of $\{\beta^{\star}\}$ in the projective space $\Proj H^1(\Sigma_g, \, {\cal O}(L_x^{-2}))\simeq\Proj^g$ under the natural projection $H^1(\Sigma_g, \, {\cal O}(L_x^{-2}))\setminus\{0\}\longrightarrow\Proj H^1(\Sigma_g, \, {\cal O}(L_x^{-2}))$ lies in the complement of the curve which is the image of the map (\ref{eqn:ext-class-projmap}). Since $g\geq 2$, the complement of a complex curve in the $g$-dimensional complex projective space $\Proj^g$ provides plenty of room for choice of deformations $E=E_t$, $t\neq 0$, of the trivial rank-two vector bundle $E_0$. Of course, $E_0$ corresponds to the trivial extension (\ref{eqn:ext-Lx-Lx-1}) or, equivalently, to the zero class $\{\beta^{\star}\}\in H^1(\Sigma_g, \, {\cal O}(L_x^{-2}))$.  \hfill $\Box$

\vspace{2ex}

 As explained before the statement of Theorem \ref{The:ES_triv-stab-def}, this result provides the final argument to the proof of the Eastwood-Singer part $(b)$ of Theorem \ref{The:E1=Einf_op-cl}. We have chosen to reproduce the approach of [ES93] in some detail because it throws up new stimulating questions of which we now mention just a few.

 The deformation behaviour of compact complex manifolds satisfying the $\partial\bar\partial$-lemma seems to be shrouded in mystery. However, deformations of twistor spaces may hold the key. Let us recall the following result of Gauduchon.

\begin{The}([Gau91, Proposition 11.b, p. 618])\label{The:Gau-twist} Every twistor space is balanced\footnote{The term {\it semi-K\"ahler} is used in [Gau91] to mean {\it balanced}.}. 

 Moreover, this is the case for any self-dual metric on the corresponding four-manifold (i.e. a balanced metric is obtained on the twistor space from any self-dual metric in a given conformal class of the base four-manifold).

\end{The}

 Combined with [ES93], this readily (re-)proves part $(a)$ of Theorem \ref{The:sG-specseq-unrel}. Indeed, by Gauduchon's Theorem \ref{The:Gau-twist}, the central fibre $Z_0$ in the Eastwood-Singer family is an sG manifold since, as a twistor space, it even has the stronger {\it balanced} property. However, $E_1(Z_0)\neq E_{\infty}(Z_0)$ as has been seen in Proposition \ref{Prop:Frolicher-nondeg}. Recall that the Iwasawa manifold had provided another example proving $(a)$ of Theorem \ref{The:sG-specseq-unrel} (cf. comments after Observation \ref{Obs:coh-Iwasawa}).

 On the other hand, Hitchin's main result in [Hit81] states that there exist only two K\"ahler twistor spaces: $\Proj^3$ and the space of flags in $\C^3$. Together with Gauduchon's theorem \ref{The:Gau-twist} and Campana's result in [Cam91b] stating that a twistor space is Moishezon if and only if it is of {\it class} ${\cal C}$, this fact will point to a relative lack of variation in the properties of twistor spaces. Thus one may legitimately ask the following

\begin{Question}\label{Question:twist-ddbar} $(a)$\, Do there exist twistor spaces that are not {\bf class} ${\cal C}$ but on which the $\partial\bar\partial$-lemma holds? 

$(b)$\, Do there exist twistor spaces whose Fr\"olicher spectral sequence degenerates at $E_1$ but on which the $\partial\bar\partial$-lemma does not hold?

$(c)$\, More generally, can one characterise the twistor spaces on which the $\partial\bar\partial$-lemma holds?

\end{Question}

 Note that if the answer to part $(a)$ of Question \ref{Question:twist-ddbar} is {\it No}, then the examples of Campana [Cam91a] and Lebrun-Poon [LP92] that proved the deformation non-openness of the {\it class} ${\cal C}$ property will also prove the deformation non-openness of the $\partial\bar\partial$-lemma property of compact complex manifolds. Indeed, in these references, families of twistor spaces are constructed in which the central fibre $Z_0$ is a Moishezon ($=$ {\it class} ${\cal C}$) twistor space (hence its algebraic dimension is maximal: $a(Z_0)=\mbox{dim}_{\C}Z_0=3$), while all the nearby fibres can be chosen to be decidedly non-Moishezon ($=$ non-{\it class} ${\cal C}$) twistor spaces, i.e. their algebraic dimension is minimal: $a(Z_t)=0$ for $t\neq 0$. Now, the Moishezon central fibre $Z_0$ satisfies the $\partial\bar\partial$-lemma, while the fibres $Z_t$ corresponding to $t\neq 0$ will not satisfy the $\partial\bar\partial$-lemma if the answer to part $(a)$ of Question \ref{Question:twist-ddbar} is {\it No}. However, independently of the eventual answer to the general question $(a)$, we can ask whether the specific fibres $Z_t$ with $t\neq 0$ in the examples mentioned above satisfy the $\partial\bar\partial$-lemma.

 In a similar vein, note that if the answer to part $(b)$ of Question \ref{Question:twist-ddbar} is {\it No}, then the Eastwood-Singer construction [ES93] that we reproduced above will also prove the deformation non-closedness of the $\partial\bar\partial$-lemma property of compact complex manifolds. We have seen that there is plenty of room for choice of stable deformations of the trivial rank-two vector bundle over $\Sigma_g$ in the Eastwood-Singer construction: the complement of a curve in the complex projective space $\Proj^g$, $g\geq 2$. It is thus natural to ask whether it is possible to make this choice in order to further ensure the validity of the $\partial\bar\partial$-lemma on every fibre $Z_t$ over $t\in\Delta^{\star}$ in the resulting family of twistor spaces.

 Note also that our Theorem \ref{The:Moi-closed} proved in [Pop10a] rules out the possibility that both parts $(a)$ and $(b)$ of Question \ref{Question:twist-ddbar} above have a {\it No} answer. Indeed, otherwise the Moishewon, {\it class} ${\cal C}$, $\partial\bar\partial$-lemma and Fr\"olicher degeneration at $E_1$ properties would be all equivalent for twistor spaces. However, if this were the case, all the fibres $Z_t$ over $t\in\Delta^{\star}$ in the Eastwood-Singer example would have to be Moishezon, hence by Theorem \ref{The:Moi-closed}, $Z_0$ would also be Moishezon. Nevertheless, $Z_0$ is not Moishezon since $E_1(Z_0)\neq E_{\infty}(Z_0)$. It is thus natural to raise part $(c)$ of Question \ref{Question:twist-ddbar}. The difficulty with answering these questions is that no condition on the base four-manifold $M$ is known that can guarantee the twistor space $Z$ of $M$ to satisfy the $\partial\bar\partial$-lemma. Does there exist an equation on $M$ similar to (\ref{eqn:H11-eqn}) that can be an indicator of the $\partial\bar\partial$-lemma property on $Z$?  

 We conclude by briefly mentioning another class of compact complex manifolds that usually provide a host of examples of various sorts: {\it nilmanifolds}. Recall that a connected, simply connected {\it real} Lie group $G$ possessing a discrete co-compact subgroup $\Gamma$ and a left-invariant complex structure $J$ defines a quotient compact complex manifold $X:=\Gamma\backslash G$ (inheriting its complex structure from $J$ by passing to the quotient) that is said to be a {\it nilmanifold} if the group $G$ is {\it nilpotent}. In the special case of a {\it complex} Lie group $G$, the quotient $X$ is a compact {\it complex parallelisable} nilmanifold (e.g. the Iwasawa manifold) of the class introduced by Wang in [Wan54], but $G$ need not be a {\it complex} Lie group in general. Putting the following two pieces of information together: 

-by [DGMS75], any compact complex manifold on which the $\partial\bar\partial$-lemma holds is {\it formal};

-by [Has89], the only {\it formal} nilmanifolds are the complex tori (i.e. those defined by an {\it abelian} Lie group $G$),

\noindent and since complex tori are clearly K\"ahler, we see that K\"ahlerness and the $\partial\bar\partial$-lemma property are equivalent conditions on nilmanifolds. Thus the class of nilmanifolds is not suitable for the study of the deformation properties of compact complex manifolds on which the $\partial\bar\partial$-lemma holds.

\vspace{3ex}

\noindent {\bf References.} \\

\noindent [AHS78]\, M. F. Atiyah, N. J. Hitchin, I. M. Singer --- {\it Self-duality in Four-Dimensional Riemannian Geometry} --- Proc. Roy. Soc. London Ser. A {\bf 362} (1978) 425-461.

\vspace{1ex}

\noindent [AB90]\, L. Alessandrini, G. Bassanelli --- {\it Small Deformations of a Class of Compact Non-K\"ahler Manifolds} --- Proc. Amer. Math. Soc. {\bf 109} (1990), no. 4, 1059–1062. 

\vspace{1ex}

\noindent [AB91a]\, L. Alessandrini, G. Bassanelli --- {\it Compact $p$-K\"ahler Manifolds} --- Geometriae Dedicata {\bf 38} (1991) 199-210.

\vspace{1ex}

\noindent [AB91b]\, L. Alessandrini, G. Bassanelli --- {\it Smooth Proper Modifications of Compact K\"ahler Manifolds} --- Proc. Internat. Workshop on Complex Analysis (Wuppertal 1990); Complex Analysis, Aspects of mathematics, {\bf E17}, Vieweg, Braunschweig (1991), 1-7.

\vspace{1ex}

\noindent [AB93]\, L. Alessandrini, G. Bassanelli --- {\it Metric Properties of Manifolds Bimeromorphic to Compact K\"ahler Spaces} --- J. Diff. Geom. {\bf 37} (1993), 95-121.

\vspace{1ex}

\noindent [AB95]\, L. Alessandrini, G. Bassanelli --- {\it Modifications of Compact Balanced Manifolds} --- C. R. Acad. Sci. Paris, t {\bf 320}, S\'erie I (1995), 1517-1522.

\vspace{1ex}

\noindent [Bes87]\, A.L. Besse --- {\it Einstein Manifolds} --- Springer, Berlin (1987).

\vspace{1ex}

\noindent [Buc99]\, N. Buchdahl --- {\it On Compact K\"ahler Surfaces} --- Ann. Inst. Fourier {\bf 49}, no. 1 (1999) 287-302.

\vspace{1ex}

\noindent [Cam91a]\, F. Campana --- {\it The Class ${\cal C}$ Is Not Stable by Small Deformations} --- Math. Ann. {\bf 290} (1991) 19-30.

\vspace{1ex}

\noindent [Cam91b]\, F. Campana --- {\it On Twistor Spaces of the Class ${\cal C}$} --- J. Diff. Geom. {\bf 33} (1991) 541-549.

\vspace{1ex}

\noindent [CE53]\, E. Calabi, B. Eckmann --- {\it A Class of Compact, Complex Manifolds Which Are Not Algebraic} --- Ann. of Math. {\bf 58} (1953) 494-500.

\vspace{1ex}

\noindent [CFGU00]\, L.A. Cordero, M. Fernandez, A. Gray, L. Ugarte --- {\it Compact Nilmanifolds with Nilpotent Complex Structures: Dolbeault Cohomology} --- Trans. Amer. Math. Soc. {\bf 352} No. 12, 5405-5433.

\vspace{1ex}

\noindent [DGMS75]\, P. Deligne, Ph. Griffiths, J. Morgan, D. Sullivan --- {\it Real Homotopy Theory of K\"ahler Manifolds} --- Invent. Math. {\bf 29} (1975), 245-274.

\vspace{1ex}

\noindent [DP04]\, J.-P. Demailly, M. Paun --- {\it Numerical Charaterization of the K\"ahler Cone of a Compact K\"ahler Manifold} --- Ann. of Math. (2) {\bf 159(3)} (2004) 1247-1274.

\vspace{1ex}

\noindent [ES93]\, M. Eastwood, M. Singer --- {\it The Fr\"ohlicher (sic) Spectral Sequence on a Twistor Space} --- J. Diff. Geom. {\bf 38} (1993) 653-669.

\vspace{1ex}

\noindent [Fuj78]\, A. Fujiki --- {\it Closedness of the Douady Spaces of Compact K\"ahler Spaces} --- Publ. RIMS, Kyoto Univ. {\bf 14} (1978), 1-52.

\vspace{1ex}

\noindent [Gau77a]\, P. Gauduchon --- {\it Le th\'eor\`eme de l'excentricit\'e nulle} --- C.R. Acad. Sc. Paris, S\'erie A, t. {\bf 285} (1977), 387-390.

\vspace{1ex}

\noindent [Gau77b]\, P. Gauduchon --- {\it Fibr\'es hermitiens \`a endomorphisme de Ricci non n\'egatif} --- Bull. Soc. Math. France {\bf 105} (1977) 113-140.

\vspace{1ex}

\noindent [Gau91]\, P. Gauduchon --- {\it Structures de Weyl et th\'eor\`emes d'annulation sur une vari\'et\'e conforme autoduale} --- Ann. Scuola Norm. Sup. Pisa Cl. Sci. (4) {\bf 18} (1991), no. 4, 563–629. 

\vspace{1ex}

\noindent [Has89]\, K. Hasegawa --- {\it Minimal Models of Nilmanifolds} --- Proc. Amer. Math. Soc. {\bf 106}, No. 1 (1989) 65-71.

\vspace{1ex}

\noindent [Hir62]\, H. Hironaka --- {\it An Example of a Non-K\"ahlerian Complex-Analytic Deformation of K\"ahlerian Complex Structures} --- Ann. of Math. (2) {\bf 75 (1)} (1962), 190-208.

\vspace{1ex}

\noindent [Hit81]\, N. J. Hitchin --- {\it K\"ahlerian Twistor Spaces} --- Proc. London Math. Soc. {\bf 43} (1981) 133-150.

\vspace{1ex}

\noindent [Hit87]\, N. J. Hitchin --- {\it The Self-duality Equations on a Riemann Surface} --- Proc. London Math. Soc. (3) {\bf 55} (1987), no.1, 59–126.

\vspace{1ex}

\noindent [Hop48]\, H. Hopf --- {\it Zur Topologie der komplexen Mannigfaltigkeiten} --- Studies and Essays Presented to R. Courant on his 60th Birthday, January 8, 1948, pp. 167–185. Interscience Publishers, Inc., New York, 1948.

\vspace{1ex}

\noindent [JS93]\, S. Ji, B. Shiffman --- {\it Properties of Compact Complex Manifolds Carrying Closed Positive Currents} --- J. Geom. Anal. {\bf 3(1)} (1993) 37-61.

\vspace{1ex}

\noindent [Kod86]\, K. Kodaira --- {\it Complex Manifolds and Deformations of Complex Structures} --- Grundlehren der Math. Wiss. {\bf 283}, Springer (1986).

\vspace{1ex}

\noindent [KS60]\, K. Kodaira, D.C. Spencer --- {\it On Deformations of Complex Analytic Structures. III. Stability Theorems for Complex Structures} --- Ann. of Math. (2) {\bf 71} No.1 (1960) 43-76. 

\vspace{1ex}

\noindent [Kur62]\, M. Kuranishi --- {\it On the Locally Complete Families of Complex Analytic Structures} --- Ann. of Math. {\bf 75}, no. 3 (1962), 536-577.

\vspace{1ex}

\noindent [Lam99]\, A. Lamari --- {\it Courants k\"ahl\'eriens et surfaces compactes} --- Ann. Inst. Fourier {\bf 49}, no. 1 (1999), 263-285.

\vspace{1ex}

\noindent [Leb86]\, C. Lebrun --- {\it On the Topology of Self-Dual $4$-Manifolds} --- Proc. Amer. Math. Soc. {\bf 98} (1986) 637-640.

\vspace{1ex}

\noindent [LP92]\, C. Lebrun, Y.-S. Poon --- {\it Twistors, K\"ahler Manifolds, and Bimeromorphic Geometry. II} --- J. Amer. Math. Soc. {\bf 5}, No. 2 (1992) 317-325.

\vspace{1ex}

\noindent [Mic83]\, M. L. Michelsohn --- {\it On the Existence of Special Metrics in Complex Geometry} --- Acta Math. {\bf 143} (1983) 261-295.

\vspace{1ex}

\noindent [Miy74]\, Y. Miyaoka --- {\it K\"ahler Metrics on Elliptic Surfaces} --- Proc. Japan Acad. {\bf 50} No. 8 (1974) 533-536.

\vspace{1ex}

\noindent [Moi67]\, B.G. Moishezon --- {\it On $n$-dimensional Compact Varieties with $n$ Algebraically Independent Meromorphic Functions} --- Amer. Math. Soc. Translations {\bf 63} (1967) 51-177.

\vspace{1ex}

\noindent [Nak75]\, I. Nakamura --- {\it Complex Parallelisable Manifolds and their Small Deformations} --- J. Diff. Geom. {\bf 10} (1975) 85-112.

\vspace{1ex}

\noindent [NR69]\, M.S. Narasimhan, S. Ramanan --- {\it Moduli of Vector Bundles on a Compact Riemann Surface} --- Ann. of Math. {\bf 89} No. 1 (1969) 14-51.

\vspace{1ex}

\noindent [NS65]\, M. S. Narasimhan, C. S. Seshadri --- {\it Stable and Unitary Vector Bundles on a Compact Riemann Surface} --- Ann. of Math. {\bf 82} (1965) 540–567.

\vspace{1ex}

\noindent [Pen76]\, R. Penrose --- {\it Nonlinear Gravitons and Curved Twistor Theory} --- General Relativity Gravitation {\bf 7} (1976) 31-52.

\vspace{1ex}

\noindent [Pop09]\, D. Popovici --- {\it Limits of Projective Manifolds under Holomorphic Deformations} --- arXiv e-print math.AG/0910.2032v1.

\vspace{1ex}

\noindent [Pop10a]\, D. Popovici --- {\it Limits of Moishezon Manifolds under Holomorphic Deformations} --- arXiv e-print math.AG/1003.3605v1.

\vspace{1ex}

\noindent [Pop10b]\, D. Popovici --- {\it Stability of Strongly Gauduchon Manifolds under Modifications} --- arXiv e-print math.CV/1009.5408v1.

\vspace{1ex}

\noindent [Sch07]\, M. Schweitzer --- {\it Autour de la cohomologie de Bott-Chern} --- arXiv e-print math. AG/0709.3528v1.

\vspace{1ex}

\noindent [Siu83]\, Y.-T. Siu --- {\it Every K3 Surface Is K\"ahler} --- Invent. Math. {\bf 73} (1983) 139-150.

\vspace{1ex}

\noindent [Tsu84]\, H. Tsuji --- {\it Complex Structures on $S^3\times S^3$} --- Tohoku Math. J (2) {\bf 36} (1984), no. 3, 351-376.

\vspace{1ex}

 \noindent [Var86]\, J. Varouchas --- {\it Sur l'image d'une vari\'et\'e k\"ahl\'erienne compacte} --- LNM {\bf 1188}, Springer (1986) 245--259.

\vspace{1ex}

\noindent [Wan54]\, H.C. Wang --- {\it Complex Parallisable (sic) Manifolds} --- Proc. Amer. Math. Soc. {\bf 5} (1954), 771-776.

\vspace{6ex}

\noindent Institut de Math\'ematiques de Toulouse, Universit\'e Paul Sabatier

\noindent 118 route de Narbonne, 31062 Toulouse Cedex 9, France

\noindent Email: popovici@math.ups-tlse.fr

\end{document}